\documentclass[11pt]{article}
\pdfoutput=1
\usepackage[utf8]{inputenc}
\usepackage[T1]{fontenc}
\usepackage{lmodern}
\usepackage{microtype}
\usepackage[letterpaper,textwidth=16cm,textheight=623.01pt,
  hcentering,top=2.5cm,headheight=12pt,headsep=25pt]{geometry}
\usepackage[shortlabels]{enumitem}

\usepackage{e-graphs-macros}
\numberwithin{equation}{section}

\usepackage{tikz-cd}
\usepackage{quiver}
\usetikzlibrary{shapes.geometric}
\tikzset{
  var/.style={inner sep=1pt},
  op/.style={draw, trapezium, trapezium angle=110, trapezium stretches=true,
             minimum width=8mm, minimum height=0pt, inner xsep=1pt, inner ysep=1.5pt,
             execute at begin node=\vphantom{$f$}},
  edge/.style={->, >=stealth},
}

\usepackage[breaklinks,hidelinks,bookmarksnumbered]{hyperref}
\usepackage[capitalize,noabbrev,nameinlink]{cleveref}

\usepackage{keytheorems}
\newkeytheorem{theorem}[within=section]
\newkeytheorem{proposition}[sibling=theorem]
\newkeytheorem{corollary}[sibling=theorem]
\newkeytheorem{lemma}[sibling=theorem]
\newkeytheorem{definition}[sibling=theorem,style=definition,qed=\qedsymbol]
\newkeytheorem{construction}[sibling=theorem,style=definition,qed=\qedsymbol,
  refname={Construction,Constructions}]
\newkeytheorem{remark}[sibling=theorem,style=remark,qed=\qedsymbol]
\newkeytheorem{example}[sibling=theorem,style=remark,qed=\qedsymbol]

\usepackage[style=alphabetic,maxbibnames=5]{biblatex}
\renewbibmacro{in:}{}

\title{How to prove equations using diagrams: \\
  A categorical and operational semantics for e-graphs}
\author{Evan Patterson}
\date{}

\begin{document}
\maketitle

\begin{abstract}
  E-graphs are a powerful algorithmic technique for automated equational
  reasoning but have resisted a conceptual understanding independent of issues
  of presentation and implementation. We propose an operational semantics for
  e-graphs based on the fundamental category-theoretic concepts of diagrams,
  limits, and initial functors. We illustrate the formalism through numerous
  example calculations and establish its soundness and completeness. Technical
  prerequisites are kept to a minimum in the main text; less standard elements
  of the mathematical theory are developed in two appendices.
\end{abstract}

\section{Introduction}
\label{sec:introduction}

Mathematical expressions and computer programs are most directly represented as
the labeled trees known as \emph{terms}. However, computer scientists have long
found reasons to encode expressions in more general mathematical and data
structures. Generalizing from trees to directed acyclic graphs (DAGs), a
\emph{term DAG} enables common subexpressions to be represented by the same
vertices, making more efficient use of space \cite{paterson1978}. Generalizing
further to equip a term DAG with an equivalence relation, an \emph{equality
  graph}, or \emph{e-graph} for short, efficiently encodes a family of terms
along with an equivalence relation on the terms
\cite{nelson1980,nelson-oppen1980}. Computing the equivalence relation's
congruence closure under term formation then leads to an algorithm to decide
quantifier-free formulas in the theory of equality with uninterpreted function
symbols. The ideas behind e-graphs, dating back to the late 1970s, received a
surge of new attention when their combination with equality saturation
\cite{tate2009} was shown to be a fast and flexible tool for equational
reasoning \cite{willsey2021}.

Standard expositions, early and recent, emphasize the data structures and
algorithms needed to implement e-graphs, such as the union-find
\cite{tarjan1975} and hash-cons data structures and congruence closure
algorithms. Data structures and algorithms are indispensable to building
practical technology based on e-graphs. Still, they do not exhaust the ways to
understand e-graphs. What might be called a ``theory of e-graphs''---a
mathematical and conceptual understanding of e-graphs, independent of any
algorithmic or implementation details---has been elusive. Besides its intrinsic
interest, a good theory can be expected to eventually lead to further practical
advances.

In this paper, we propose a mathematical semantics for e-graphs. The semantics
is, in the first place, a \emph{categorical semantics}, meaning that
category-theoretic concepts---categories, functors, and so forth---are used to
give mathematical content to syntactical objects like terms or, in this case,
term graphs. Characteristic of categorical semantics is its \emph{invariance} to
presentation: not depending on how the objects in question are presented or
implemented. In this spirit, we deliberately avoid proposing a single definition
of ``e-graph.'' Rather, we describe a means of representing equations using
category-theoretic diagrams that is conceptually compatible with how e-graphs
encode equations, but makes no commitments about presentation.

In addition to being categorical, our semantics is primarily \emph{operational}.
That might seem surprising since categorical semantics are often exclusively
denotational, beginning with Lawvere's reconstruction of algebraic theories as
categories with finite products \cite{lawvere1963,crole1993,adamek2011}. Indeed,
we will use algebraic theories qua categories in \cref{sec:cartesian-setting}.
Yet our guiding principle throughout is to analyze not what e-graphs are as
\emph{objects} but what their valid transformations are as \emph{morphisms}.
This is an operational semantics inasmuch as it gives mathematical content to
operations performed on e-graphs. The operational aspect, like the denotational
one, is invariant: the set of allowed operations is defined ``all at once'' and
without reference to a presentation, rather than being generated by a minimal
set of atomic steps as in a small-step operational semantics. Our operational
semantics also does not aim at reduction to a final value (or set of values).
This is in contrast to many operational semantics for programming languages, but
it is in keeping with software systems like \texttt{egg} \cite{willsey2021} and
\texttt{egglog} \cite{zhang2023} that are not fixed procedures but allow
flexible schedules of rules to be run against an e-graph.

The thesis of this paper is that a categorical and operational semantics for
e-graphs can be found in the concept of an \emph{initial functor}. Consider the
term graphs below, which for the moment we take to be labeled graphs.
\begin{equation*}
  \begin{tikzcd}[column sep=small]
    & x & \\
    {y_1} && {y_2} \\
    w && z
    \arrow["f"', from=1-2, to=2-1]
    \arrow["f", from=1-2, to=2-3]
    \arrow["g"', from=2-1, to=3-1]
    \arrow["h", from=2-3, to=3-3]
  \end{tikzcd}
  \qquad\to\qquad
  \begin{tikzcd}[column sep=small]
    & x & \\
    & y \\
    w && z
    \arrow["f"', from=1-2, to=2-2]
    \arrow["g"', from=2-2, to=3-1]
    \arrow["h", from=2-2, to=3-3]
  \end{tikzcd}
\end{equation*}
Both term graphs encode the same two terms, $g(f(x))$ and $h(f(x))$. The map
between them, collapsing both edges labeled $f$ on the left to the single edge
labeled $f$ on the right, is a homomorphism of labeled graphs. It is an example
of a \emph{congruence reduction}, which is, in turn, one type of valid
transformation of term graphs, where ``valid'' means not changing the terms or
equations implicit in the graph. It is then natural to ask: what are \emph{all}
the valid transformations? To our knowledge, an elegant answer to this question
cannot be given using only the language of graphs and graph homomorphisms. If,
instead, we see the map as a functor (between categories freely generated by
graphs), we can recognize an example of a known class of functors, the
\emph{initial} functors. As we will explain, initial functors are exactly those
that preserve the system of equations encoded by a diagram; in
category-theoretic language, restricting the shape of a diagram along an initial
functor preserves the diagram's limit, whenever it exists.

The operational semantics for e-graphs thus rests on an analysis of when a
functor is guaranteed to maintain, in a suitable sense, the \emph{universal
  property} defining a limit. Universal properties are an essential---arguably
\emph{the} essential---idea in category theory. It is for this reason that
category theory is the natural language for a theory of e-graphs.

\paragraph{Outline}

\cref{sec:representation} lays the groundwork by explaining how to represent
systems of equations using category-theoretic diagrams. \cref{sec:reasoning}
justifies this representation by introducing morphisms between diagrams and
singling out the initial ones as those preserving the systems of equations.
Rules for, and examples of, calculating with initial morphisms are given. So
\cref{sec:representation,sec:reasoning} concern representation and reasoning,
respectively. They are restricted to categories with no extra structure, which
in logical terms correspond to algebraic theories with only unary operations.
This setting is very limited practically, but it suffices to illustrate the main
ideas. \cref{sec:cartesian-setting} extends the formalism to categories with
finite products, corresponding to the classical setting of algebraic theories
(with multiary operations). Finally, \cref{sec:outlook} offers a brief outlook
on future directions. Two technical appendices develop the background
mathematics.

\paragraph{Background}

We assume familiarity with the basic concepts of category theory, particularly
categories, functors, diagrams, and limits. Further concepts not covered by most
textbooks, namely initial functors and comprehensive factorization, are reviewed
in \cref{sec:initial-functors}. As for computer science, readers familiar with
term graphs or e-graphs will find connections with those concepts, but no prior
knowledge of them is needed to make sense of this paper.

\paragraph{Related work}

Two groups have recently proposed category-theoretic accounts of e-graphs
\cite{biondo2025,tiurin2025}. Without going into specifics, both define an
``e-hypergraph'' to be a labeled hypergraph equipped with extra structure that,
directly or indirectly, specifies an equivalence relation. Neither
representation allows taking initial functors to be morphisms of e-graphs,
the main idea and novelty of this paper.

Motivated by ``Tonti diagrams'' in computational physics, the author and
collaborators have previously used category-theoretic diagrams to represent
systems of equations---here, \emph{differential} equations---and initial
functors to formalize weak equivalence of diagrams
\cite{patterson2023,carlson2024}. The underlying mathematics dates back to the
early 1970s, with key results due to Par{\'e} \cite{pare1973} and Street--Walters
\cite{street1973}. The novelty of this work is to recognize the connection between
initial functors and the operational aspect of e-graphs.

\section{Representing equations using diagrams}
\label{sec:representation}

A system of equations can be represented by a \emph{diagram}. Here, we use the
word ``diagram'' in the precise sense of category theory. Recall that a
\define{diagram} in a category $\cat{C}$ consists of a small category $\cat{J}$,
called the diagram's \define{shape}, along with a functor
$D: \cat{J} \to \cat{C}$.

While diagrams can encode equations in different ways, in this paper a diagram
specifies a system of equations through its \emph{limit}. The limit may or may
not exist in the target category $\cat{C}$, but the limit always exists when the
diagram is interpreted by sets and functions via a functor $F: \cat{C} \to \Set$.
A simple first example illustrates these points.

\begin{example}[An equalizer diagram]
  Consider a square in some category $\cat{C}$:
  \begin{equation}
    \label{eq:square-diagram}
    \begin{tikzcd}
      {x: X} & {y: Y} \\
      {w: W} & {z:Z}
      \arrow["f", from=1-1, to=1-2]
      \arrow["h"', from=1-1, to=2-1]
      \arrow["g", from=1-2, to=2-2]
      \arrow["k"', from=2-1, to=2-2]
    \end{tikzcd}.
  \end{equation}
  First of all, what diagram $D: \cat{J} \to \cat{C}$ is meant by this picture? In
  general, we use lowercase letters for objects in the shape $\cat{J}$ and
  uppercase letters for objects in the target category
  $\cat{C}$.\footnote{Explicitly naming the objects in a diagram's shape is
    unconventional among category theorists, but doing so makes space for the
    \emph{variables} that show up in term graphs.} So $w,x,y,z$ are objects in
  $\cat{J}$ and $W,X,Y,Z$ are objects in $\cat{C}$. As for morphisms, we leave
  the morphisms in the shape nameless and use lowercase letters for morphisms in
  the target category. So, in this example, $\cat{C}$ must contain morphisms
  $f: X \to Y$, $g: Y \to Z$, $h: X \to W$, and $k: W \to Z$.

  In the absence of further assumptions, the diagram is not completely
  determined by its depiction \eqref{eq:square-diagram}. The freest, or most
  general, assumption would take both $\cat{J}$ and $\cat{C}$ to be categories
  freely generated by graphs. The diagram would then be given explicitly as the
  functor
  \begin{equation*}
    \cat{J} = \left\{
      \begin{tikzcd}[cramped, sep=small]
        x & y \\
        w & z
        \arrow[from=1-1, to=1-2]
        \arrow[from=1-1, to=2-1]
        \arrow[from=1-2, to=2-2]
        \arrow[from=2-1, to=2-2]
      \end{tikzcd}
    \right\}
    \xto{\ D\ }
    \left\{
      \begin{tikzcd}[cramped, sep=small]
        X & Y \\
        W & Z
        \arrow["f", from=1-1, to=1-2]
        \arrow["h"', from=1-1, to=2-1]
        \arrow["g", from=1-2, to=2-2]
        \arrow["k"', from=2-1, to=2-2]
      \end{tikzcd}
    \right\} = \cat{C}.
  \end{equation*}

  Whatever the choice of target category $\cat{C}$, the diagram represents the
  system of equations
  \begin{equation} \label{eq:equalizer-system-1}
    \left\{
      x: X,\, y: Y,\, w: W,\, z: Z
      \;\middle|\;
      f(x) = y,\, g(y) = z,\, h(x) = w,\, k(w) = z
    \right\}.
  \end{equation}
  This system of equations is equivalently but more compactly written as
  \begin{equation} \label{eq:equalizer-system-2}
    \left\{
      x: X
      \;\middle|\;
      g(f(x)) = k(h(x))
    \right\},
  \end{equation}
  although that is not the literal reading of the diagram.

  The formal meaning of the ``system of equations represented by a diagram''
  comes from the diagram's limit, which characterizes the system's solutions by
  a universal property. The limit can be taken ``semantically'' or
  ``syntactically.'' Working semantically, for any functor
  $F: \cat{C} \to \Set$, the composite functor $F \circ D: \cat{J} \to \Set$ is
  a diagram of sets and functions whose limit is (in canonical bijection with)
  both \eqref{eq:equalizer-system-1} and \eqref{eq:equalizer-system-2}, now
  interpreted as actual sets.\footnote{This reading of
    \eqref{eq:equalizer-system-1} and \eqref{eq:equalizer-system-2} abuses
    notation by conflating the abstract object (``type'') $X$ with the set
    $F(X)$ and the morphism (``operation'') $f$ with the function $F(f)$.} Or,
  staying in the realm of syntax, we can try to evaluate the limit directly in
  $\cat{C}$. If $\cat{C}$ is the free category on a directed square, then the
  limit does not exist, reflecting that the diagram's equations do not hold in
  $\cat{C}$. On the other hand, if $\cat{C}$ is finitely presented as
  \begin{equation*}
    \cat{C} \coloneqq \left\langle
      \begin{tikzcd}[cramped, sep=small]
        X & Y \\
        W & Z
        \arrow["f", from=1-1, to=1-2]
        \arrow["h"', from=1-1, to=2-1]
        \arrow["g", from=1-2, to=2-2]
        \arrow["k"', from=2-1, to=2-2]
      \end{tikzcd}
      \;\Biggm|\;
      g \circ f = k \circ h
      \right\rangle,
  \end{equation*}
  then the limit of the diagram $D$ exists and its apex is just $X$, reflecting
  that the diagram's equations are now true in $\cat{C}$. Later,
  we will see how to turn this observation into a technique for proving
  equations in a finitely presented category.
\end{example}

Systematizing the example, a diagram $D: \cat{J} \to \cat{C}$ is \define{free}
if its shape is a free category, i.e., $\cat{J} = \Path(G)$ is freely generated
by a graph $G$.\footnote{When speaking precisely of a \define{graph}, we always
  mean the category theorist's notion---what a graph theorist might call a
  ``directed multigraph with self-loops.'' In \cref{sec:cartesian-setting} we
  use ``multigraph'' to mean something else.} Note that a diagram being free
implies nothing about its target category. The examples in this section and the
next will be exclusively of free diagrams. We emphasize, however, that nothing
in the theory requires this; in future work, we may consider diagram shapes
containing nontrivial isomorphisms as a more faithful model of the equivalence
relation in an e-graph.

By the free-forgetful adjunction between graphs and categories, a free diagram
in a category $\cat{C}$ is nothing other than a graph labeled by objects and
morphisms in $\cat{C}$. The important connection to category theory is not
through the data of the diagram itself, but through the notions of limit and, as
explained in \cref{sec:reasoning}, of limit-preserving changes of diagram shape.

\subsection{Other ways to represent equations with diagrams}

It should be noted that this way of representing equations with diagrams subtly
differs from others with which the reader may be more familiar.

\begin{remark}[Comparison with commutative diagrams]
  \label{rem:commutative-diagrams}
  Category theorists often use commutative diagrams to assert equations between
  morphisms. While it is not always carefully defined, this phrase does have a
  precise meaning: a diagram $D: \cat{J} \to \cat{C}$ \define{commutes}, or is
  \define{commutative}, if for every pair of morphisms in $\cat{J}$ with the
  same domain and codomain, their images in $\cat{C}$ are equal.

  The senses in which limits and commutativity represent equations through a
  diagram are \emph{not} the same. As a small counterexample, consider the free
  diagrams
  \begin{equation*}
    \begin{tikzcd}[column sep=small]
      & {x: X} & \\
      {y_1: Y} && {y_2: Y} \\
      & {z: Z}
      \arrow["f"', from=1-2, to=2-1]
      \arrow["f", from=1-2, to=2-3]
      \arrow["g"', from=2-1, to=3-2]
      \arrow["h", from=2-3, to=3-2]
    \end{tikzcd}
    \qquad\text{and}\qquad
    \begin{tikzcd}
      {x: X} \\
      {y: Y} \\
      {z: Z}
      \arrow["f"', from=1-1, to=2-1]
      \arrow["g"', curve={height=12pt}, from=2-1, to=3-1]
      \arrow["h", curve={height=-12pt}, from=2-1, to=3-1]
    \end{tikzcd}
  \end{equation*}
  in some category where $g \circ f = h \circ f$ but $g \neq h$. Then both diagrams have
  the same limit, namely $X$, but the first diagram commutes while the second
  does not.
\end{remark}

\begin{remark}[Comparison with cyclic term graphs]
  \label{rem:cyclic-term-graphs}
  The graph rewriting community generalizes the notion of term graph from a DAG
  to an arbitrary graph, regarding it as a system of recursion equations
  \cite[e.g.,][\mbox{Figure 6}]{ariola1996}. Cycles are unwound into infinite
  terms with finitely many distinct subterms, called \emph{rational terms}
  \cite{corradini1998}. Our formalism interprets diagrams differently, with no
  special treatment of cycles. Each morphism in the diagram's shape introduces
  an equational constraint between a pair of variables, where in the case of an
  endomorphism the two variables happen to be the same (see
  \cref{ex:monoid-presentation} below).
\end{remark}

\section{Reasoning with initial morphisms between diagrams}
\label{sec:reasoning}

What makes the foregoing representation useful, compared to innumerable other
ways of encoding terms and equations in graph-like structures, is that its
notion of \emph{morphism} leads naturally to a way of reasoning about equations.

A \define{morphism of diagrams} in a common category $\cat{C}$, say from
$D_0: \cat{J} \to \cat{C}$ to $D: \cat{K} \to \cat{C}$, is a functor
$R: \cat{J} \to \cat{K}$ making the triangle
\begin{equation*}
  \begin{tikzcd}[column sep=small]
    {\cat{J}} && {\cat{K}} \\
    & {\cat{C}}
    \arrow["R", from=1-1, to=1-3]
    \arrow["{D_0}"', from=1-1, to=2-2]
    \arrow["D", from=1-3, to=2-2]
  \end{tikzcd}
\end{equation*}
commute.\footnote{A more flexible definition of morphism between diagrams
  permits the triangle to be filled by a natural transformation, instead of
  strictly commuting. However, it turns out \cite{carlson2024} that this
  ostensibly more general starting point leads to the same relation of weak
  equivalence on diagrams, so we need not consider it here.} In other words, the
\define{category of diagrams} in $\cat{C}$ is the slice category
$\cat{Cat}/\cat{C}$.

An arbitrary morphism of diagrams need not preserve the systems of equations
represented by the diagrams, but it does give a canonical comparison in one
direction. Assuming the limits in question exist, the restricted family of
projections $\pi_{Rj}: \lim D \to D(R(j))$, for $j \in \cat{J}$, forms a cone over
$D_0 = D \circ R$ and thus, by the universal property of the limit of $D_0$, there
exists a unique morphism
\begin{equation} \label{eq:limit-comparison}
  \textstyle
  \lim_{\cat{K}} D \to \lim_{\cat{J}} D_0 = \lim_{\cat{J}} (D \circ R)
\end{equation}
factoring through the projections of $\lim D_0$. In general, this comparison map
is not invertible.

A functor $R: \cat{J} \to \cat{K}$ is \define{initial} if, for \emph{every}
diagram $D: \cat{K} \to \cat{C}$ in an arbitrary category $\cat{C}$, the canonical
comparison map \eqref{eq:limit-comparison} is an isomorphism, whenever the
limits exist (either limit existing if the other does). For other, equivalent
ways to define an initial functor, see \cref{def:initial-functor} in
\cref{sec:initial-functors}. We denote an initial functor with a double-headed
arrow, $R: \cat{J} \initialto \cat{K}$, and we extend the terminology to say
that a morphism of diagrams is \define{initial} if its underlying functor is
initial.

Intuitively, the initial morphisms of diagrams are those that preserve the
system of equations represented by the domain or codomain diagrams, as the
following example illustrates.

\begin{example}[A congruence reduction]
  \label{ex:congruence-reduction}
  Resuming the example from \cref{rem:commutative-diagrams}, the morphism of
  diagrams
  \begin{equation*}
    \begin{tikzcd}[column sep=small]
      & {x: X} & \\
      {y_1: Y} && {y_2: Y} \\
      & {z: Z}
      \arrow["f"', from=1-2, to=2-1]
      \arrow["f", from=1-2, to=2-3]
      \arrow["g"', from=2-1, to=3-2]
      \arrow["h", from=2-3, to=3-2]
    \end{tikzcd}
    \qquad\longinitialto\qquad
    \begin{tikzcd}
      {x: X} \\
      {y: Y} \\
      {z: Z}
      \arrow["f"', from=1-1, to=2-1]
      \arrow["g"', curve={height=12pt}, from=2-1, to=3-1]
      \arrow["h", curve={height=-12pt}, from=2-1, to=3-1]
    \end{tikzcd},
  \end{equation*}
  sending both $y_1$ and $y_2$ in the domain's shape to $y$ in the codomain's
  shape, is initial. To see this, notice that the limit of both diagrams is the
  equalizer of $g \circ f$ and $h \circ f$, whenever it exists. The intuitive
  justification is that the systems of equations presented,
  \begin{equation*}
    \left\{
      x: X,\, y_1: Y,\, y_2: Y,\, z: Z
      \;\middle|\;
      f(x) = y_1,\, f(x) = y_2,\, g(y_1) = z,\, h(y_2) = z
    \right\}
  \end{equation*}
  and
  \begin{equation*}
    \left\{
      x: X,\, y: Y,\, z: Z
      \;\middle|\;
      f(x) = y,\, g(y) = z,\, h(y) = z
    \right\},
  \end{equation*}
  are equivalent to each other, and also to $\{x: X \mid g(f(x)) = h(f(x))\}$.
\end{example}

\subsection{Weak equivalence of diagrams}
\label{sec:weak-equivalence}

Initial functors, and thus also initial morphisms of diagrams, are closed under
composition and identities (\cref{prop:initial-functor-closure-properties}), so
two diagrams being ``related by an initial functor'' is a transitive and reflexive
relation. It is not symmetric: although every isomorphism of
categories---indeed, every equivalence of categories
(\cref{prop:initial-functor-closure-properties})---is an initial functor, there
are plenty of initial functors that are not invertible (see
\cref{ex:congruence-reduction}). Those are the interesting ones. Yet the
relation on diagrams of ``representing equivalent systems of equations'' should
plainly be an equivalence relation; said differently, we want to treat initial
morphisms \emph{as if} they were invertible. This problem is solved by formally
inverting the initial morphisms.

A \define{weak equivalence} of diagrams will be another name for an initial
morphism of diagrams. Two diagrams $D$ and $D'$ in some category are
\define{weakly equivalent}, denoted $D \sim D'$, if they can be connected by a
zig-zag of weak equivalences.

We say \emph{weak} equivalence because such morphisms are not equivalences in
any familiar 2-category of diagrams. The phrase ``weak equivalence'' is also
deliberately suggestive of homotopy theory. While we will not make serious use
of homotopy-theoretic techniques, we note in passing one characteristic property
\cite[\S{I.2}]{gabriel1967}: \emph{a priori}, weakly equivalent diagrams are
connected only by a zig-zag of arbitrary finite length, but in fact, a cospan
always suffices. To restate this:

\begin{proposition}[Weak equivalence, cospan form]
  Two diagrams $D: \cat{J} \to \cat{C}$ and $D': \cat{J}' \to \cat{C}$ are weakly
  equivalent if and only if there exists a diagram $P: \cat{K} \to \cat{C}$ and
  initial functors $R: \cat{J} \initialto \cat{K}$ and
  $R': \cat{J}' \initialto \cat{K}$ making the triangles commute:
  \begin{equation*}
    \begin{tikzcd}
      {\cat{J}} & {\cat{K}} & {\cat{J}'} \\
      & {\cat{C}}
      \arrow["R", two heads, from=1-1, to=1-2]
      \arrow["D"', from=1-1, to=2-2]
      \arrow["P"{description}, from=1-2, to=2-2]
      \arrow["{R'}"', two heads, from=1-3, to=1-2]
      \arrow["{D'}", from=1-3, to=2-2]
    \end{tikzcd}.
  \end{equation*}
\end{proposition}
\begin{proof}
  Follows from the stability of initial functors, and hence of initial morphisms
  of diagrams, under pushouts (\cref{prop:initial-functor-closure-properties}).
\end{proof}

We record this fact only to emphasize that weak equivalence can be defined ``all
at once,'' without recourse to a finite but unbounded iteration. Directly
calculating the intermediate diagram $P: \cat{K} \to \cat{C}$ does not seem to
be computationally feasible.

\subsection{Rules of inference}
\label{sec:rules}

Although initial functors and weak equivalence are invariant notions applicable
to diagrams of any size, when calculating it is convenient to have a library of
motifs that can be applied mechanically, without having to check initiality at each
step. In this section, we record the basic motifs and explain what it means to
apply them.

In the motifs that follow, all diagrams are in a fixed but arbitrary category
$\cat{C}$.
\begin{itemize}
  \item \define{Operation introduction}: for any morphism $f: X \to Y$ in $\cat{C}$,
    the diagram morphism
    \begin{equation*}
      x: X
      \qquad\longinitialto\qquad
      x: X \xto{f} y:Y
    \end{equation*}
    is initial.
  \item \define{Congruence reduction}: for any morphism $f: X \to Y$ in $\cat{C}$,
    the diagram morphism
    \begin{equation*}
      \begin{tikzcd}[column sep=small]
        & {x:X} & \\
        {y_1: Y} && {y_2: Y}
        \arrow["f"', from=1-2, to=2-1]
        \arrow["f", from=1-2, to=2-3]
      \end{tikzcd}
      \qquad\longinitialto\qquad
      \begin{tikzcd}
        {x:X} \\
        {y:Y}
        \arrow["f"', from=1-1, to=2-1]
      \end{tikzcd},
    \end{equation*}
    sending both $y_1$ and $y_2$ to $y$, is initial.
  \item \define{Equation elimination}: for any morphism $f: X \to Y$ in
    $\cat{C}$, the diagram morphism
    \begin{equation*}
      \begin{tikzcd}
        {x: X} \\
        {y:Y}
        \arrow["f", curve={height=-12pt}, from=1-1, to=2-1]
        \arrow["f"', curve={height=12pt}, from=1-1, to=2-1]
      \end{tikzcd}
      \qquad\longinitialto\qquad
      \begin{tikzcd}
        {x:X} \\
        {y:Y}
        \arrow["f"', from=1-1, to=2-1]
      \end{tikzcd}
    \end{equation*}
    is initial.
  \item \define{Factorization expansion}: for any morphisms $f: X \to Y$ and
    $g: Y \to Z$ in $\cat{C}$, the diagram morphism
    \begin{equation*}
      x: X \xto{g \circ f} z: Z
      \qquad\longinitialto\qquad
      x: X \xto{f} y: Y \xto{g} z: Z
    \end{equation*}
    is initial. Similarly, for any object $X \in \cat{C}$, the diagram morphism
    \begin{equation*}
      x_1: X \xto{\id_X} x_2: X
      \qquad\longinitialto\qquad
      x: X,
    \end{equation*}
    sending both $x_1$ and $x_2$ to $x$, is initial.
\end{itemize}
All of these functors are initial since the limit of every diagram above is
simply $X$.

Each motif gives rise to a rule of inference that can be applied to transform
one diagram into another that is weakly equivalent. We will not go so far as to
propose a formal system for such derivations, as our focus is rather on
mathematical semantics. Still, the justification for our calculations is that
initial functors, and thus also initial morphisms of diagrams, are stable under
pushout with respect to arbitrary morphisms
(\cref{prop:initial-functor-closure-properties}). An example illustrates how to
apply rules by computing pushouts.

\begin{example}[Applying congruence reduction]
  The initial morphism of diagrams asserted directly in
  \cref{ex:congruence-reduction} can be constructed by applying the congruence
  reduction rule:
  \newsavebox{\ruleDomain}
  \begin{lrbox}{\ruleDomain}
    \begin{tikzcd}[cramped, column sep=small]
      & {x:X} & \\
      {y_1:Y} && {y_2:Y}
      \arrow["f"', from=1-2, to=2-1]
      \arrow["f", from=1-2, to=2-3]
    \end{tikzcd}
  \end{lrbox}
  \newsavebox{\ruleCodomain}
  \begin{lrbox}{\ruleCodomain}
    \begin{tikzcd}[cramped]
      {x:X} \\
      {y:Y}
      \arrow["f"', from=1-1, to=2-1]
    \end{tikzcd}
  \end{lrbox}
  \newsavebox{\exDomain}
  \begin{lrbox}{\exDomain}
    \begin{tikzcd}[cramped, column sep=small]
      & {x:X} & \\
      {y_1:Y} && {y_2:Y} \\
      & {z:Z}
      \arrow["f"', from=1-2, to=2-1]
      \arrow["f", from=1-2, to=2-3]
      \arrow["g"', from=2-1, to=3-2]
      \arrow["h", from=2-3, to=3-2]
    \end{tikzcd}
  \end{lrbox}
  \newsavebox{\exCodomain}
  \begin{lrbox}{\exCodomain}
    \begin{tikzcd}[cramped]
      {x:X} \\
      {y:Y} \\
      {z:Z}
      \arrow["f"', from=1-1, to=2-1]
      \arrow["g"', curve={height=12pt}, from=2-1, to=3-1]
      \arrow["h", curve={height=-12pt}, from=2-1, to=3-1]
    \end{tikzcd}
  \end{lrbox}
  \begin{equation*}
    \begin{tikzcd}[column sep=huge, row sep=large]
      {\left\{\usebox{\ruleDomain}\right\}} & {\left\{\usebox{\ruleCodomain}\right\}} \\
      {\left\{\usebox{\exDomain}\right\}} & {\left\{\usebox{\exCodomain}\right\}}
      \arrow[two heads, from=1-1, to=1-2]
      \arrow[from=1-1, to=2-1]
      \arrow[from=1-2, to=2-2]
      \arrow[two heads, from=2-1, to=2-2]
      \arrow["\lrcorner"{anchor=center, pos=0.125, rotate=180}, draw=none, from=2-2, to=1-1]
    \end{tikzcd}
    \qedhere
  \end{equation*}
\end{example}

To show that two diagrams are weakly equivalent, we are allowed to apply these
rules in either direction. In doing so we can derive new rules. Particularly
important is \define{equation introduction} for equations between morphisms in
the target category $\cat{C}$. As a representative case, suppose $g \circ f = k \circ h$
in $\cat{C}$. Applying the rules above, there is a zig-zag of weak equivalences
\begin{equation*}
  \begin{tikzcd}
    {x: X} \\
    {}
  \end{tikzcd}
  \quad\longinitialto\quad
  \begin{tikzcd}
    {x: X} \\
    {z: Z}
    \arrow["{g \circ f = k \circ h}"{description}, from=1-1, to=2-1]
  \end{tikzcd}
  \quad\longinitialfrom\quad
  \begin{tikzcd}
    {x: X} \\
    {z: Z}
    \arrow["{g \circ f}"', curve={height=12pt}, from=1-1, to=2-1]
    \arrow["{k \circ h}", curve={height=-12pt}, from=1-1, to=2-1]
  \end{tikzcd}
  \quad\longinitialto\quad
  \begin{tikzcd}[cramped, sep=small]
    & {x: X} & \\
    {y: Y} && {w: W} \\
    & {z: Z}
    \arrow["f"', from=1-2, to=2-1]
    \arrow["h", from=1-2, to=2-3]
    \arrow["g"', from=2-1, to=3-2]
    \arrow["k", from=2-3, to=3-2]
  \end{tikzcd}.
\end{equation*}
Hence, whenever $g \circ f = k \circ h$ in $\cat{C}$, the two diagrams are
weakly equivalent:
\begin{equation*}
  {x: X}
  \qquad\sim\qquad
  \begin{tikzcd}[cramped, sep=small]
    & {x: X} & \\
    {y: Y} && {w: W} \\
    & {z: Z}
    \arrow["f"', from=1-2, to=2-1]
    \arrow["h", from=1-2, to=2-3]
    \arrow["g"', from=2-1, to=3-2]
    \arrow["k", from=2-3, to=3-2]
  \end{tikzcd}.
\end{equation*}

It remains to explain how to \emph{begin} a derivation, assuming nothing of the
target category. Suppose we start with a \define{singleton} diagram, whose shape
is the terminal category, say $\{x: X\}$. The limit of this diagram exists
trivially, being $X$ itself. Any weakly equivalent diagram will continue to have
$X$ as its limit. In particular, any parallel paths in such a diagram that
originate at the object corresponding to $x$ must commute. We use this technique
in the calculations that follow. It will become more powerful when we extend to
categories with finite products (\cref{sec:cartesian-setting}), since we can
then start with any finite discrete diagram, the shape's generating objects
acting as ``free variables.''

\subsection{Example calculations}

We show by example how to prove equations using the calculus of weak
equivalence.

\begin{example}[Pasting commutative squares]
  The pasting law for commutative squares says that two commutative squares with
  a common edge can be juxtaposed to produce a commutative rectangle. The
  notation makes this seem obviously true---that's what makes it good
  notation---but it still needs a proof. The calculus of weak equivalence makes
  the proof completely mechanical.

  Suppose that the two squares below commute:
  \begin{equation*}
    \begin{tikzcd}
      X & Y \\
      {X'} & {Y'}
      \arrow["f", from=1-1, to=1-2]
      \arrow["h"', from=1-1, to=2-1]
      \arrow["k", from=1-2, to=2-2]
      \arrow["{f'}"', from=2-1, to=2-2]
    \end{tikzcd}
    \qquad\qquad
    \begin{tikzcd}
      Y & Z \\
      {Y'} & {Z'}
      \arrow["g", from=1-1, to=1-2]
      \arrow["k"', from=1-1, to=2-1]
      \arrow["\ell", from=1-2, to=2-2]
      \arrow["{g'}"', from=2-1, to=2-2]
    \end{tikzcd}.
  \end{equation*}
  Formally, take the target category $\cat{C}$ to be generated by six objects,
  called $X,X',Y,Y',Z,Z'$, and seven morphisms, called $f,f',g,g',h,k,\ell$,
  subject to those two equations.

  Begin with the singleton diagram $\{x: X\}$, representing the domain of the
  equation to be proved. Recognize that we can introduce the first equation at
  $x$:
  \begin{equation*}
    \begin{tikzcd}
      {x: X} & {y: Y} \\
      {x': X'} & {y': Y'}
      \arrow["f", from=1-1, to=1-2]
      \arrow["h"', from=1-1, to=2-1]
      \arrow["k", from=1-2, to=2-2]
      \arrow["{f'}"', from=2-1, to=2-2]
    \end{tikzcd}.
  \end{equation*}
  Then recognize that we can introduce the second equation at $y$:
  \begin{equation*}
    \begin{tikzcd}
      {x: X} & {y: Y} & {z:Z} & \\
      {x': X'} & {y': Y'} & {\tilde y': Y'} & {z': Z'}
      \arrow["f", from=1-1, to=1-2]
      \arrow["h"', from=1-1, to=2-1]
      \arrow["g", from=1-2, to=1-3]
      \arrow["k", from=1-2, to=2-2]
      \arrow["k", from=1-2, to=2-3]
      \arrow["\ell", from=1-3, to=2-4]
      \arrow["{f'}"', from=2-1, to=2-2]
      \arrow["{g'}"', from=2-3, to=2-4]
    \end{tikzcd}.
  \end{equation*}
  Finally, apply congruence reduction at $y$:
  \begin{equation*}
    \begin{tikzcd}
      {x: X} & {y: Y} & {z: Z} \\
      {x': X'} & {y': Y'} & {z': Z'}
      \arrow["f", from=1-1, to=1-2]
      \arrow["h"', from=1-1, to=2-1]
      \arrow["g", from=1-2, to=1-3]
      \arrow["k", from=1-2, to=2-2]
      \arrow["\ell", from=1-3, to=2-3]
      \arrow["{f'}"', from=2-1, to=2-2]
      \arrow["{g'}"', from=2-2, to=2-3]
    \end{tikzcd}.
  \end{equation*}
  Since the limit of the original diagram is $X$, the limit of the final, weakly
  equivalent diagram remains $X$. Consequently the outer paths around the
  rectangle commute, which was the equation in $\cat{C}$ to be proved:
  \begin{equation*}
    \begin{tikzcd}
      X & Y & Z \\
      {X'} & {Y'} & {Z'}
      \arrow["f", from=1-1, to=1-2]
      \arrow["h"', from=1-1, to=2-1]
      \arrow["g", from=1-2, to=1-3]
      \arrow["\ell", from=1-3, to=2-3]
      \arrow["{f'}"', from=2-1, to=2-2]
      \arrow["{g'}"', from=2-2, to=2-3]
    \end{tikzcd}.
    \qedhere
  \end{equation*}
\end{example}

The next example is adapted from Nelson and Oppen's original work introducing
e-graphs \cite[\S{6}]{nelson1980}, \cite{nelson-oppen1980}.

\begin{example}[A monoid presentation]
  \label{ex:monoid-presentation}
  Consider a monoid with an element $f$ satisfying the equations
  \begin{equation*}
    f^3 = \id = f^5.
  \end{equation*}
  That is, take the target category $\cat{C}$ to be generated by a single
  object, say $\bullet$, and an endomorphism, $f: \bullet \to \bullet$, satisfying those two
  equations.

  Begin with the singleton diagram $\{x_0: \bullet\}$ and introduce both defining
  equations, noting that since each equation has an identity on one side, they
  appear as cycles:
  \begin{equation*}
    \begin{tikzcd}
      {x_2} & {x_0} & {x_4'} & {x_3'} \\
      {x_1} & {x_1'} & {x_2'}
      \arrow[from=1-1, to=1-2]
      \arrow[from=1-2, to=2-1]
      \arrow[from=1-2, to=2-2]
      \arrow[from=1-3, to=1-2]
      \arrow[from=1-4, to=1-3]
      \arrow[from=2-1, to=1-1]
      \arrow[from=2-2, to=2-3]
      \arrow[from=2-3, to=1-4]
    \end{tikzcd}.
  \end{equation*}
  We omit the types of the objects and the labels of the arrows, these always
  being the object $\bullet$ and the generating morphism $f$, respectively. We now
  compute the \emph{congruence closure}, iteratively applying congruence
  reduction until that is no longer possible. The first two reductions, at $x_0$
  and then at $x_1$, are:
  \begin{equation*}
    \longinitialto\quad
    \begin{tikzcd}
      {x_2} & {x_0} & {x_4'} & {x_3'} \\
      & {x_1} & {x_2'}
      \arrow[from=1-1, to=1-2]
      \arrow[from=1-2, to=2-2]
      \arrow[from=1-3, to=1-2]
      \arrow[from=1-4, to=1-3]
      \arrow[from=2-2, to=1-1]
      \arrow[from=2-2, to=2-3]
      \arrow[from=2-3, to=1-4]
    \end{tikzcd}
    \quad\longinitialto\quad
    \begin{tikzcd}
      {x_0} & {x_4'} & {x_3'} \\
      {x_1} & {x_2}
      \arrow[from=1-1, to=2-1]
      \arrow[from=1-2, to=1-1]
      \arrow[from=1-3, to=1-2]
      \arrow[from=2-1, to=2-2]
      \arrow[from=2-2, to=1-1]
      \arrow[from=2-2, to=1-3]
    \end{tikzcd}.
  \end{equation*}
  Applying congruence reduction five more times, at $x_2$, $x_0$, $x_1$, $x_0$,
  and $x_0$, gives the sequence:
  \begin{equation*}
    \longinitialto\quad
    \begin{tikzcd}
      {x_0} & {x_4'} \\
      {x_1} & {x_2}
      \arrow[curve={height=-6pt}, from=1-1, to=1-2]
      \arrow[from=1-1, to=2-1]
      \arrow[curve={height=-6pt}, from=1-2, to=1-1]
      \arrow[from=2-1, to=2-2]
      \arrow[from=2-2, to=1-1]
    \end{tikzcd}
    \quad\longinitialto\quad
    \begin{tikzcd}
      {x_0} \\
      {x_1} & {x_2}
      \arrow[from=1-1, to=2-1]
      \arrow[curve={height=-12pt}, from=2-1, to=1-1]
      \arrow[from=2-1, to=2-2]
      \arrow[from=2-2, to=1-1]
    \end{tikzcd}
    \quad\longinitialto\quad
    \begin{tikzcd}
      {x_0} \\
      {x_1}
      \arrow[from=1-1, to=1-1, loop, in=60, out=120, distance=5mm]
      \arrow[curve={height=-6pt}, from=1-1, to=2-1]
      \arrow[curve={height=-6pt}, from=2-1, to=1-1]
    \end{tikzcd}
    \quad\longinitialto\quad
    \begin{tikzcd}
      {x_0}
      \arrow[from=1-1, to=1-1, loop, in=60, out=120, distance=5mm]
      \arrow[from=1-1, to=1-1, loop, in=240, out=300, distance=5mm]
    \end{tikzcd}
    \quad\longinitialto\quad
    \begin{tikzcd}
      {x_0}
      \arrow[from=1-1, to=1-1, loop, in=60, out=120, distance=5mm]
    \end{tikzcd}.
  \end{equation*}
  The final diagram is the self-loop over $f$. We have proved the equation
  \begin{equation*}
    f = \id,
  \end{equation*}
  since the limit of the final diagram is that of the original one, namely the
  unique object $\bullet$. In other words, the monoid presented is the trivial monoid.
\end{example}

\subsection{Soundness and completeness}
\label{sec:soundness-completeness}

The calculus of weak equivalence is put on a firmer foundation by showing that
it is sound and complete for reasoning about systems of equations under
\emph{any} interpretation by sets and functions. In logical jargon, one might
say that it is sound and complete for equational reasoning with ``uninterpreted
function symbols.'' Soundness holds essentially by the definition of an initial
functor (\cref{def:initial-functor}). Completeness is more involved. Here we
state only the result, deferring the proofs to \cref{sec:initial-functors} and
references therein.

\begin{theorem}
  Two diagrams $D: \cat{J} \to \cat{C}$ and $D': \cat{J}' \to \cat{C}$ are
  weakly equivalent if and only if, for every set-valued functor
  $F: \cat{C} \to \Set$, the limits $\lim_{\cat{J}} (F \circ D)$ and
  $\lim_{\cat{J}'} (F \circ D')$ are in bijection, naturally in $F$.
\end{theorem}
\begin{proof}
  These are two of the four equivalent conditions in
  \cref{thm:image-copresheaf-isomorphism}, a more refined result.
\end{proof}

We caution that we have not proved the completeness of the specific rules in
\cref{sec:rules}, though we expect some such result should hold when precisely
formulated.

\section{The cartesian setting: universal algebra via diagrams}
\label{sec:cartesian-setting}

The categorical and operational semantics for e-graphs so far exists in the
simplest possible setting, that of categories absent any further properties or
structure. In the dictionary of categorical logic, a bare category corresponds
to an algebraic theory with only \emph{unary} function symbols. This is rather
restrictive. The classical setting for equational reasoning---universal
algebra---is that of algebraic theories, whose categorical counterparts are
categories with finite products. Historically, universal algebra was the
starting point for both categorical logic \cite{lawvere1963} and e-graphs
\cite{nelson1980}.

In this section, we extend the operational semantics from bare categories to
categories with finite products. Surprisingly, this turns out to be, in a sense,
trivial: in the definitions and results used in
\cref{sec:representation,sec:reasoning}, it is enough to insert the adjective
``cartesian'' wherever it makes sense. Full justification of this claim is
deferred to \cref{sec:cartesian-initial-functors}, where initial functors,
comprehensive factorization, and limits are shown to interact well with finite
products. Here we focus on the operational aspect of calculating with diagrams
valued in an algebraic theory.

\subsection{Cartesian diagrams and their morphisms}

For brevity, we say that a category is \define{cartesian} if it has finite
products and that a functor is \define{cartesian} if it preserves finite
products. A diagram $D: \cat{J} \to \cat{C}$ in a cartesian category $\cat{C}$
is \define{cartesian} if its shape $\cat{J}$ is a cartesian category and the
functor $D$ is itself cartesian.

Cartesian categories and functors are familiar; cartesian diagrams less so, as
limits or colimits are not usually indexed by categories assumed to have
products. When taking limits in a cartesian category, there is no harm in making
the diagrams cartesian too, because a cone over a cartesian diagram is
uniquely determined at a product of objects by its legs at the factors
(\cref{lem:cones-cartesian-diagram}). Moreover, any of the diagrams from
\cref{sec:representation,sec:reasoning} can be transported to the cartesian
setting without altering its limit (\cref{prop:limit-free-product-completion}).
But cartesian diagrams are more expressive, allowing diagram shapes to specify
multiary operations and diagram morphisms to assign arbitrary terms.

Cartesian diagrams, like ordinary ones, are conveniently generated by
combinatorial data. The appropriate data is a \define{multigraph}:\footnote{This
  terminology is chosen to make the analogy ``graphs : categories :: multigraphs :
  multicategories'' hold.} a set of \define{vertices} and a set of
\define{multiedges}, along with an assignment to each multiedge of a \emph{list} of
vertices, called its \define{source}, and a single vertex, its \define{target}.
In more loaded jargon, a multigraph is a \define{signature}, its vertices are
\define{sorts}, and its multiedges are \define{function symbols}. A
\define{morphism of multigraphs} consists of a pair of functions between the
vertex sets and multiedge sets, preserving the source and target assignments.
Rearranging the data slightly, the category of multigraphs is seen to be a
presheaf topos \cite[\S{1}]{walters1989}.

We say that a cartesian diagram is \define{free} if its shape is a cartesian
category freely generated by a multigraph. There are several ways to explicitly
construct the free cartesian category on a multigraph $M$
\cite[e.g.,][\S{2}]{walters1989}, but it is enough to know that a free cartesian
diagram $D: F(M) \to \cat{C}$ corresponds\footnote{The nature of this
  correspondence, a 2-categorical universal property, is the subject of Walters'
  paper \cite{walters1989}.} to a multigraph morphism $M \to U(\cat{C})$ into the
underlying multigraph of $\cat{C}$. The latter is nothing other than a
\define{term graph}, taking that phrase to mean a labeled multigraph. Moreover,
the limit of the free cartesian diagram is completely determined by the term
graph (\cref{prop:cones-free-cartesian-diagram}). As before, all examples of
cartesian diagrams in this paper will be free, but freeness is not an assumption
of the theory.

The target category may be any cartesian category but is often given by finite
presentation. We need not be too precise about how to present a cartesian
category. Roughly, a presentation consists of a signature and a set of equations
between terms in context over the signature. Adopting the perspective of
categorical logic \cite{lawvere1963,crole1993,adamek2011}, we identify a small
cartesian category with an algebraic theory, where the objects and morphisms of
the category correspond to the types and operations of the theory, respectively.

A first example is adapted from Nelson and Oppen \cite[\S{6}]{nelson1980},
\cite{nelson-oppen1980}; along with \cref{ex:monoid-presentation}, this covers
the two toy examples from the original work on e-graphs.

\begin{example}
  \label{ex:term-graph}
  Let $\cat{C}$ be the signature with (i.e., the free cartesian category
  generated by) a type $X$ and a binary operation ${f: (X, X) \to X}$. Then the
  terms $f(a, b)$ and $f(f(a, b), b)$ with free variables $a: X$ and $b: X$, as
  well as the equation $f(a, b) = a$, are encoded by the term graph (i.e., the
  free cartesian diagram) in $\cat{C}$:
  \begin{equation*}
    \begin{tikzpicture}[baseline=(current bounding box.center)]
      \node[var] (a) at (-1, 0) {$a: X$};
      \node[var] (b) at (1, 0) {$b: X$};
      \node[op] (f1) at (0, -0.9) {$f$};
      \node[op] (f2) at (0, -2.1) {$f$};
      \node[var] (c) at (0, -2.9) {$c: X$};
      \draw[edge] (a) -- ([xshift=-2mm]f1.north);
      \draw[edge] (b) -- ([xshift=2mm]f1.north);
      \draw[edge] (f1.south) to[out=-90, in=-70, looseness=1.2] (a);
      \draw[edge] (a) to[out=-105, in=90] ([xshift=-2mm]f2.north);
      \draw[edge] (b) to[out=-75, in=90] ([xshift=2mm]f2.north);
      \draw[edge] (f2.south) -- (c);
    \end{tikzpicture}.
  \end{equation*}
  Later, we will see how to simplify this diagram using congruence closure.
\end{example}

The notions of morphism from \cref{sec:reasoning} extend straightforwardly to
the cartesian setting. A \define{morphism of cartesian diagrams} is a morphism of
diagrams whose underlying functor is cartesian. In other words, the
\define{category of cartesian diagrams} in a cartesian category $\cat{C}$ is the
slice category $\cat{Cart}/\cat{C}$. A \define{cartesian initial functor} is nothing
other than a functor that is both cartesian and initial. Similarly, a morphism
of cartesian diagrams is \define{initial} if its underlying functor is initial.
Such morphisms are also called \define{weak equivalences}. Two cartesian
diagrams are \define{weakly equivalent} if they can be connected by a zig-zag of
weak equivalences.

\subsection{Rules of inference}

The motifs and rules introduced in \cref{sec:rules} become more general in the
cartesian setting, effectively allowing a context comprising a single free
variable to be replaced by an arbitrary context, i.e., a finite list of
variables.

In the motifs that follow, let $\cat{C}$ be a cartesian category.
\begin{itemize}
  \item \define{Operation introduction}: for any morphism
    $f: (X_1, \dots, X_k) \to Y$ in $\cat{C}$, the morphism of cartesian
    diagrams
    \begin{equation*}
      \begin{tikzpicture}[baseline=(current bounding box.center)]
        \node[var] (x1) at (-1.25, 0) {$x_1: X_1$};
        \node (xmid) at (0, 0) {$\dots$};
        \node[var] (xk) at (1.25, 0) {$x_k: X_k$};
      \end{tikzpicture}
      \qquad\longinitialto\qquad
      \begin{tikzpicture}[baseline=(current bounding box.center)]
        \node[var] (x1) at (-1.25, 0) {$x_1: X_1$};
        \node (xmid) at (0, 0) {$\dots$};
        \node[var] (xk) at (1.25, 0) {$x_k: X_k$};
        \node[var] (y) at (0, -2) {$y: Y$};
        \node[op] (f) at (0, -1) {$f$};
        \draw[edge] (x1) -- ([xshift=-2mm]f.north);
        \draw[edge] (xmid) -- (f.north);
        \draw[edge] (xk) -- ([xshift=2mm]f.north);
        \draw[edge] (f.south) -- (y);
      \end{tikzpicture}
    \end{equation*}
    is initial.
  \item \define{Congruence reduction}: for any morphism
    $f: (X_1, \dots, X_k) \to Y$ in $\cat{C}$, the morphism of cartesian
    diagrams,
    \begin{equation*}
      \begin{tikzpicture}[baseline=(current bounding box.center)]
        \node[var] (x1) at (-1.25, 0) {$x_1: X_1$};
        \node (xmid) at (0, 0) {$\dots$};
        \node[var] (xk) at (1.25, 0) {$x_k: X_k$};
        \node[op] (f1) at (-0.6, -1) {$f$};
        \node[op] (f2) at (0.6, -1) {$f$};
        \node[var] (y1) at (-0.6, -2) {$y_1: Y$};
        \node[var] (y2) at (0.6, -2) {$y_2: Y$};
        \draw[edge] (x1) -- ([xshift=-2mm]f1.north);
        \draw[edge] (xmid) -- (f1.north);
        \draw[edge] (xk) -- ([xshift=2mm]f1.north);
        \draw[edge] (x1) -- ([xshift=-2mm]f2.north);
        \draw[edge] (xmid) -- (f2.north);
        \draw[edge] (xk) -- ([xshift=2mm]f2.north);
        \draw[edge] (f1.south) -- (y1);
        \draw[edge] (f2.south) -- (y2);
      \end{tikzpicture}
      \qquad\longinitialto\qquad
      \begin{tikzpicture}[baseline=(current bounding box.center)]
        \node[var] (x1) at (-1.25, 0) {$x_1: X_1$};
        \node (xmid) at (0, 0) {$\dots$};
        \node[var] (xk) at (1.25, 0) {$x_k: X_k$};
        \node[var] (y) at (0, -2) {$y: Y$};
        \node[op] (f) at (0, -1) {$f$};
        \draw[edge] (x1) -- ([xshift=-2mm]f.north);
        \draw[edge] (xmid) -- (f.north);
        \draw[edge] (xk) -- ([xshift=2mm]f.north);
        \draw[edge] (f.south) -- (y);
      \end{tikzpicture},
    \end{equation*}
    sending both $y_1$ and $y_2$ on the left-hand side to $y$ on the right-hand
    side, is initial.
  \item \define{Equation elimination}: for any morphism
    $f: (X_1, \dots, X_k) \to Y$ in $\cat{C}$, the morphism of cartesian
    diagrams
    \begin{equation*}
      \begin{tikzpicture}[baseline=(current bounding box.center)]
        \node[var] (x1) at (-1.25, 0) {$x_1: X_1$};
        \node (xmid) at (0, 0) {$\dots$};
        \node[var] (xk) at (1.25, 0) {$x_k: X_k$};
        \node[op] (f1) at (-0.6, -1) {$f$};
        \node[op] (f2) at (0.6, -1) {$f$};
        \node[var] (y) at (0, -2) {$y: Y$};
        \draw[edge] (x1) -- ([xshift=-2mm]f1.north);
        \draw[edge] (xmid) -- (f1.north);
        \draw[edge] (xk) -- ([xshift=2mm]f1.north);
        \draw[edge] (x1) -- ([xshift=-2mm]f2.north);
        \draw[edge] (xmid) -- (f2.north);
        \draw[edge] (xk) -- ([xshift=2mm]f2.north);
        \draw[edge] (f1.south) -- (y);
        \draw[edge] (f2.south) -- (y);
      \end{tikzpicture}
      \qquad\longinitialto\qquad
      \begin{tikzpicture}[baseline=(current bounding box.center)]
        \node[var] (x1) at (-1.25, 0) {$x_1: X_1$};
        \node (xmid) at (0, 0) {$\dots$};
        \node[var] (xk) at (1.25, 0) {$x_k: X_k$};
        \node[var] (y) at (0, -2) {$y: Y$};
        \node[op] (f) at (0, -1) {$f$};
        \draw[edge] (x1) -- ([xshift=-2mm]f.north);
        \draw[edge] (xmid) -- (f.north);
        \draw[edge] (xk) -- ([xshift=2mm]f.north);
        \draw[edge] (f.south) -- (y);
      \end{tikzpicture}
    \end{equation*}
    is initial.
  \item \define{Factorization expansion}: for any factorization of a morphism in
    $\cat{C}$, the morphism from the cartesian diagram indexed by a single
    morphism to the cartesian diagram indexed by the shape of the factorization
    is initial. For example, given morphisms
    \begin{equation*}
      f: (X_1, X_2) \to Y_1,\quad
      g: (X_3, X_4) \to Y_2, \quad
      h: (Y_1, Y_2) \to Z
    \end{equation*}
    in $\cat{C}$, the morphism of cartesian diagrams
    \begin{equation*}
      \begin{tikzpicture}[baseline=(current bounding box.center)]
        \node[var] (x1) at (-2.25, 0) {$x_1: X_1$};
        \node[var] (x2) at (-0.75, 0) {$x_2: X_2$};
        \node[var] (x3) at (0.75, 0) {$x_3: X_3$};
        \node[var] (x4) at (2.25, 0) {$x_4: X_4$};
        \node[op] (hfg) at (0, -1.25) {$h \circ (f \times g)$};
        \node[var] (z) at (0, -2.5) {$z: Z$};
        \draw[edge] (x1) -- ([xshift=-6mm]hfg.north);
        \draw[edge] (x2) -- ([xshift=-2mm]hfg.north);
        \draw[edge] (x3) -- ([xshift=2mm]hfg.north);
        \draw[edge] (x4) -- ([xshift=6mm]hfg.north);
        \draw[edge] (hfg.south) -- (z);
      \end{tikzpicture}
      \quad\longinitialto\quad
      \begin{tikzpicture}[baseline=(current bounding box.center)]
        \node[var] (x1) at (-2.25, 0) {$x_1: X_1$};
        \node[var] (x2) at (-0.75, 0) {$x_2: X_2$};
        \node[var] (x3) at (0.75, 0) {$x_3: X_3$};
        \node[var] (x4) at (2.25, 0) {$x_4: X_4$};
        \node[op] (f) at (-1.5, -1) {$f$};
        \node[op] (g) at (1.5, -1) {$g$};
        \node[var] (y1) at (-1.5, -2) {$y_1: Y_1$};
        \node[var] (y2) at (1.5, -2) {$y_2: Y_2$};
        \node[op] (h) at (0, -3) {$h$};
        \node[var] (z) at (0, -4) {$z: Z$};
        \draw[edge] (x1) -- ([xshift=-2mm]f.north);
        \draw[edge] (x2) -- ([xshift=2mm]f.north);
        \draw[edge] (x3) -- ([xshift=-2mm]g.north);
        \draw[edge] (x4) -- ([xshift=2mm]g.north);
        \draw[edge] (f.south) -- (y1);
        \draw[edge] (g.south) -- (y2);
        \draw[edge] (y1) -- ([xshift=-2mm]h.north);
        \draw[edge] (y2) -- ([xshift=2mm]h.north);
        \draw[edge] (h.south) -- (z);
      \end{tikzpicture}
    \end{equation*}
    is initial.
\end{itemize}
To see that the functors above are initial, notice that restricting along each
one preserves the limit of any compatibly shaped cartesian diagram.\footnote{By
  \cref{prop:cartesian-initial-functor}, it suffices to consider
  \emph{cartesian} diagrams in the definition of initial functor.} Indeed, by
\cref{prop:cones-free-cartesian-diagram}, each limit is simply the product
${X_1 \times \cdots \times X_k}$ (where, in the final example, $k=4$).

The motifs give rise to rules, applied by taking pushouts. A technical hitch
compared to \cref{sec:rules} is that (1-categorical) pushouts of cartesian
categories generally do not exist. The simplest fix is to take pushouts in the
category of small categories with \emph{chosen} finite products and strictly
product-preserving functors. Within this cocomplete category, initial functors
are stable under pushout, a corollary of
\cref{thm:cartesian-comprehensive-factorization}.

\begin{example}[Applying congruence reduction]
  Returning to \cref{ex:term-graph}, the morphism of cartesian diagrams,
  \begin{equation*}
    \begin{tikzpicture}[baseline=(current bounding box.center)]
      \node[var] (a) at (-1, 0) {$a: X$};
      \node[var] (b) at (1, 0) {$b: X$};
      \node[op] (f1) at (0, -0.9) {$f$};
      \node[op] (f2) at (0, -2.1) {$f$};
      \node[var] (c) at (0, -2.9) {$c: X$};
      \draw[edge] (a) -- ([xshift=-2mm]f1.north);
      \draw[edge] (b) -- ([xshift=2mm]f1.north);
      \draw[edge] (f1.south) to[out=-90, in=-70, looseness=1.2] (a);
      \draw[edge] (a) to[out=-105, in=90] ([xshift=-2mm]f2.north);
      \draw[edge] (b) to[out=-75, in=90] ([xshift=2mm]f2.north);
      \draw[edge] (f2.south) -- (c);
    \end{tikzpicture}
    \qquad\longinitialto\qquad
    \begin{tikzpicture}[baseline=(current bounding box.center)]
      \node[var] (a) at (-1, 0) {$a: X$};
      \node[var] (b) at (1, 0) {$b: X$};
      \node[op] (f) at (0, -0.9) {$f$};
      \draw[edge] (a) -- ([xshift=-2mm]f.north);
      \draw[edge] (b) -- ([xshift=2mm]f.north);
      \draw[edge] (f.south) to[out=-90, in=-70, looseness=1.2] (a);
    \end{tikzpicture},
  \end{equation*}
  sending both $a$ and $c$ on the left-hand side to $a$ on the right-hand side,
  is initial. It is an application of congruence reduction. The diagram on the
  right cannot be reduced further.
\end{example}

We can again derive new rules, valid in the sense of preserving weak
equivalence, by applying the motifs in either direction. Notably, a generalized
\define{equation introduction} rule can be derived by the same argument as in
\cref{sec:rules}. Suppose, for example, that morphisms
\begin{equation*}
  f: (X_1, X_2) \to Y, \qquad
  g: Y \to Z, \qquad
  h: (X_1, X_2) \to Z
\end{equation*}
in a cartesian category $\cat{C}$ satisfy $h = g \circ f$. Then the rule says
that there is a weak equivalence of cartesian diagrams in $\cat{C}$:
\begin{equation*}
  \begin{tikzpicture}[baseline=(current bounding box.center)]
    \node[var] (x1) at (-0.8, 0) {$x_1: X_1$};
    \node[var] (x2) at (0.8, 0) {$x_2: X_2$};
  \end{tikzpicture}
  \qquad\sim\qquad
  \begin{tikzpicture}[baseline=(current bounding box.center)]
    \node[var] (x1) at (-0.8, 0) {$x_1: X_1$};
    \node[var] (x2) at (0.8, 0) {$x_2: X_2$};
    \node[op] (h) at (-0.8, -1.95) {$h$};
    \node[op] (f) at (0.8, -1) {$f$};
    \node[var] (y) at (0.8, -2) {$y: Y$};
    \node[op] (g) at (0.8, -2.9) {$g$};
    \node[var] (z) at (0, -3.9) {$z: Z$};
    \draw[edge] (x1) -- ([xshift=-2mm]h.north);
    \draw[edge] (x2) -- ([xshift=2mm]h.north);
    \draw[edge] (x1) -- ([xshift=-2mm]f.north);
    \draw[edge] (x2) -- ([xshift=2mm]f.north);
    \draw[edge] (f.south) -- (y);
    \draw[edge] (y) -- (g.north);
    \draw[edge] (g.south) -- (z);
    \draw[edge] (h.south) -- (z);
  \end{tikzpicture}.
\end{equation*}

Finally, the assumption-free starting point for a derivation generalizes from a
singleton diagram to a discrete diagram. Recall that a multigraph is
\define{discrete} if its set of multiedges is empty. A cartesian diagram is
\define{discrete} if its shape is freely generated by a discrete multigraph. The
limit of a cartesian diagram that is finitely generated and discrete, say
$\{x_1 : X_1, \dots, x_k : X_k\}$, always exists, being the product
$X_1 \times \cdots \times X_k$. Any weakly equivalent cartesian diagram will have the same
limit, from which equations can be deduced. The next section demonstrates this
technique.

\subsection{Example calculation: a property of well-behaved lenses}

Let $\cat{C}$ be the algebraic theory of a \define{well-behaved} lens. That is,
$\cat{C}$ is the algebraic theory with signature comprising
\begin{itemize}[noitemsep]
  \item sorts $S$ (``source'') and $V$ (``view''),
  \item operations $\lensget: S \to V$ and $\lensput: (S, V) \to S$,
\end{itemize}
and two axioms:
\begin{itemize}[noitemsep]
  \item \textsc{PutGet}: $s: S,\, v: V \vdash \lensget(\lensput(s, v)) = v$.
  \item \textsc{GetPut}: $s: S \vdash \lensput(s, \lensget(s)) = s$.
\end{itemize}
A \define{very well-behaved} or \define{lawful} lens satisfies the further
axiom:
\begin{itemize}
  \item \textsc{PutPut}:
    $s: S,\, v,v': V \vdash \lensput(\lensput(s, v), v') = \lensput(s, v')$.
\end{itemize}
The \textsc{PutPut} law is independent of the other two. However, a special
case of the \textsc{PutPut} equation, with $v = v'$, \emph{is} derivable from
the first two laws:
\begin{itemize}
  \item \textsc{PutTwice}:
    $s: S,\, v: V \vdash \lensput(\lensput(s, v), v) = \lensput(s, v)$.
\end{itemize}
We are going to prove this by diagrammatic reasoning. The derivation is fairly
trivial, though the fact is noted in the lens literature \cite{fischer2015}.

Begin with the discrete cartesian diagram $\{s: S,\, v: V\}$. Since the
right-hand side of the target equation is $\lensput(s, v)$, a natural first step
is to introduce the $\lensput$ operation:
\begin{equation*}
  \begin{tikzpicture}[baseline=(current bounding box.center)]
    \node[var] (s) at (-0.6, 0) {$s: S$};
    \node[var] (v) at (0.6, 0) {$v: V$};
  \end{tikzpicture}
  \qquad\longinitialto\qquad
  \begin{tikzpicture}[baseline=(current bounding box.center)]
    \node[var] (s) at (-0.6, 0) {$s: S$};
    \node[var] (v) at (0.6, 0) {$v: V$};
    \node[op] (put) at (0, -1) {$\lensput$};
    \node[var] (s2) at (0, -2) {$s': S$};
    \draw[edge] (s) -- ([xshift=-2mm]put.north);
    \draw[edge] (v) -- ([xshift=2mm]put.north);
    \draw[edge] (put.south) -- (s2);
  \end{tikzpicture}.
\end{equation*}
Recognizing a subterm of the \textsc{PutGet} law, we introduce that equation at
$s$ and $v$, then apply congruence reduction of $\lensput$ at $s$ and $v$:
\begin{equation*}
  \begin{tikzpicture}[baseline=(current bounding box.center)]
    \node[var] (s) at (-0.6, 0) {$s: S$};
    \node[var] (v) at (0.6, 0) {$v: V$};
    \node[op] (put) at (0, -1) {$\lensput$};
    \node[var] (s2) at (0, -2) {$s': S$};
    \draw[edge] (s) -- ([xshift=-2mm]put.north);
    \draw[edge] (v) -- ([xshift=2mm]put.north);
    \draw[edge] (put.south) -- (s2);
  \end{tikzpicture}
  \quad\sim\quad
  \begin{tikzpicture}[baseline=(current bounding box.center)]
    \node[var] (s) at (-0.6, 0) {$s: S$};
    \node[var] (v) at (0.6, 0) {$v: V$};
    \node[op] (put1) at (-0.6, -1) {$\lensput$};
    \node[op] (put2) at (0.6, -1) {$\lensput$};
    \node[var] (s2) at (-0.6, -2) {$s': S$};
    \node[var] (s3) at (0.6, -2) {$s'': S$};
    \node[op] (get) at (0.6, -2.9) {$\lensget$};
    \draw[edge] (s) -- ([xshift=-2mm]put1.north);
    \draw[edge] (v) -- ([xshift=2mm]put1.north);
    \draw[edge] (s) -- ([xshift=-2mm]put2.north);
    \draw[edge] (v) -- ([xshift=2mm]put2.north);
    \draw[edge] (put1.south) -- (s2);
    \draw[edge] (put2.south) -- (s3);
    \draw[edge] (s3) -- (get.north);
    \draw[edge, overlay] (get.south) to[out=-60, in=-10, looseness=1.6] (v.east);
    \path (1.95, -3.7); 
  \end{tikzpicture}
  \quad\longinitialto\quad
  \begin{tikzpicture}[baseline=(current bounding box.center)]
    \node[var] (s) at (-0.6, 0) {$s: S$};
    \node[var] (v) at (0.6, 0) {$v: V$};
    \node[op] (put) at (0, -1) {$\lensput$};
    \node[var] (s2) at (0, -2) {$s': S$};
    \node[op] (get) at (0, -2.9) {$\lensget$};
    \draw[edge] (s) -- ([xshift=-2mm]put.north);
    \draw[edge] (v) -- ([xshift=2mm]put.north);
    \draw[edge] (put.south) -- (s2);
    \draw[edge] (s2) -- (get.north);
    \draw[edge, overlay] (get.south) to[out=-60, in=-20, looseness=1.0] (v.east);
    \path (1.5, -3.5); 
  \end{tikzpicture}.
\end{equation*}
(We could have saved a few steps by immediately introducing the \textsc{PutGet}
equation, but we wish to emphasize the mechanizability of such reasoning.) Next,
we recognize a subterm of the \textsc{GetPut} law (replacing $s$ with $s'$), so
we introduce that equation at $s'$:
\begin{equation*}
  \begin{tikzpicture}[baseline=(current bounding box.center)]
    \node[var] (s) at (-0.6, 0) {$s: S$};
    \node[var] (v) at (0.6, 0) {$v: V$};
    \node[op] (put) at (0, -1) {$\lensput$};
    \node[var] (s2) at (0, -2) {$s': S$};
    \node[op] (get) at (0, -2.9) {$\lensget$};
    \draw[edge] (s) -- ([xshift=-2mm]put.north);
    \draw[edge] (v) -- ([xshift=2mm]put.north);
    \draw[edge] (put.south) -- (s2);
    \draw[edge] (s2) -- (get.north);
    \draw[edge, overlay] (get.south) to[out=-60, in=-20, looseness=1.0] (v.east);
    \path (1.5, -3.5); 
  \end{tikzpicture}
  \quad\sim\quad
  \begin{tikzpicture}[baseline=(current bounding box.center)]
    \node[var] (s) at (-0.6, 0) {$s: S$};
    \node[var] (v) at (0.6, 0) {$v: V$};
    \node[op] (put) at (0, -1) {$\lensput$};
    \node[var] (s2) at (0, -2) {$s': S$};
    \node[op] (get1) at (1.55, -2.9) {$\lensget$};
    \node[op] (get2) at (0.55, -2.9) {$\lensget$};
    \node[var] (v2) at (0.55, -3.8) {$v': V$};
    \node[op] (put2) at (0, -4.6) {$\lensput$};
    \draw[edge] (s) -- ([xshift=-2mm]put.north);
    \draw[edge] (v) -- ([xshift=2mm]put.north);
    \draw[edge] (put.south) -- (s2);
    \draw[edge] (s2) -- (get1.north);
    \draw[edge] (s2) -- (get2.north);
    \draw[edge] (get2.south) -- (v2);
    \draw[edge] (s2) -- ([xshift=-2mm]put2.north);
    \draw[edge] (v2) -- ([xshift=2mm]put2.north);
    \draw[edge, overlay] (get1.south) to[out=-45, in=-10, looseness=1.2] (v.east);
    \draw[edge, overlay] (put2.south) to[out=-120, in=180, looseness=1.2] (s2.west);
    \path (2.7, -3.5) (-1.0, -5.1); 
  \end{tikzpicture}
  \quad\longinitialto\quad
  \begin{tikzpicture}[baseline=(current bounding box.center)]
    \node[var] (s) at (-0.6, 0) {$s: S$};
    \node[var] (v) at (0.6, 0) {$v: V$};
    \node[op] (put) at (0, -1) {$\lensput$};
    \node[var] (s2) at (0, -2) {$s': S$};
    \node[op] (get) at (1.2, -2.9) {$\lensget$};
    \node[op] (put2) at (0, -3.9) {$\lensput$};
    \draw[edge] (s) -- ([xshift=-2mm]put.north);
    \draw[edge] (v) -- ([xshift=2mm]put.north);
    \draw[edge] (put.south) -- (s2);
    \draw[edge] (s2) -- (get.north);
    \draw[edge] (s2) -- ([xshift=-2mm]put2.north);
    \draw[edge] (v) .. controls (0.8, -1.6) and (0.7, -2.9) .. ([xshift=2mm]put2.north);
    \draw[edge, overlay] (get.south) to[out=-45, in=-10, looseness=1.2] (v.east);
    \draw[edge, overlay] (put2.south) to[out=-120, in=180, looseness=1.2] (s2.west);
    \path (2.35, -3.5) (-1.0, -4.4); 
  \end{tikzpicture}.
\end{equation*}
The final morphism is again a congruence reduction, this time of $\lensget$ at
$s'$.

The resulting cartesian diagram contains the target equation, though it takes a
moment to see that. Since the first and last cartesian diagrams are weakly
equivalent, and since the limit of the first diagram is the product ${S \times V}$,
the limit of the last diagram is also ${S \times V}$. Consequently, in the last
diagram, any two parallel morphisms in the shape starting at $s \times v$
evaluate to equal morphisms in $\cat{C}$. Taking the endpoint $s'$ and following
the cycle at $s'$ (for fixed $v$) either zero times or once, we deduce the
\textsc{PutTwice} law:
\begin{equation*}
  s: S,\, v: V \vdash \lensput(s, v) = \lensput(\lensput(s, v), v).
\end{equation*}

\section{Outlook}
\label{sec:outlook}

We have seen how the concept of initial functor leads to a method of proving
equations by manipulating diagrams. We regard this method as a categorical and
operational semantics for e-graphs, as demonstrated by reproducing calculations
from the e-graphs literature.

While we have so far treated only the classical setting of e-graphs, a good
categorical semantics should suggest generalizations of e-graphs to new or
nonstandard settings. We mention two possible directions. First, it seems likely
that the extension of initial functors and comprehensive factorization from mere
categories to cartesian categories (\cref{sec:cartesian-initial-functors}) will
extend further to cartesian equipments or, what is closely related, cartesian
fibrations. In the dictionary of categorical logic, these categorical structures
correspond to enhanced notions of algebraic theory having first-class functions
\emph{and} relations. We thus expect to find a categorical and operational
semantics for the recent combination of e-graphs and Datalog in \texttt{egglog}
\cite{zhang2023}. Another promising direction is extending to \emph{essentially
  algebraic} or \emph{generalized algebraic} theories, which correspond to
categories with all finite limits. As far as we know, e-graph software has not
yet been developed for this setting.

\paragraph{Acknowledgments}

I acknowledge support from the Air Force Office of Scientific Research (AFOSR)
Young Investigator Program (YIP) through Award FA9550-23-1-0133. I thank my
colleague Tim Hosgood for carefully reading and commenting on the series of blog
posts that preceded this paper.

\printbibliography[heading=bibintoc]

\appendix
\crefalias{section}{appendix}

\section{Review of initial functors and comprehensive factorization}
\label{sec:initial-functors}

A functor $F: \cat{C} \to \cat{D}$ \define{preserves limits} when, for every
diagram $D: \cat{J} \to \cat{C}$, the canonical comparison map
\begin{equation*}
  \textstyle F(\lim_{\cat{J}} D) \to \lim_{\cat{J}} (F \circ D),
\end{equation*}
furnished by the universal property of $\lim(F \circ D)$, is an isomorphism,
assuming the limits in question exist \cite[\mbox{Exercise V.4.5}]{maclane1998}.
That term is standard, yet there is another, dual sense in which a functor could
be said to ``preserve limits'' that varies the \emph{shape} of a diagram rather
than its target. A functor that does not change the limit of a diagram when the
diagram is restricted along it is called \emph{initial}. In this appendix, which
is purely expository, we review initial functors and the factorization system
into which they fit. Here and throughout the paper, all categories are assumed
to be locally small and shapes of diagrams to be small.

\subsection{Initial functors: definition and properties}

Initial functors are characterized by the following definition/theorem.

\begin{definition}[Initial functor]
  \label{def:initial-functor}
  A functor $R: \cat{J} \to \cat{K}$ between small categories is \define{initial},
  here denoted $R: \cat{J} \initialto \cat{K}$, if any of the following
  equivalent conditions holds:
  \begin{enumerate}[(i)]
    \item \emph{Restriction preserves limits}: for every diagram
      $D: \cat{K} \to \cat{C}$ in an arbitrary category $\cat{C}$, the canonical
      comparison map
      \begin{equation*}
        \textstyle \lim_{\cat{K}} D \to \lim_{\cat{J}} (D \circ R),
      \end{equation*}
      given by the universal property of $\lim(D \circ R)$, is an isomorphism,
      each limit existing if the other does;
    \item \emph{Restriction preserves limits of sets}: condition (i) holds for
      $\cat{C} = \Set$, i.e., for every diagram of sets $D: \cat{K} \to \Set$,
      the canonical comparison $\lim_{\cat{K}} D \to \lim_{\cat{J}} (D \circ R)$
      is a bijection;
    \item \emph{Commas are connected}: for every object $k \in \cat{K}$, the comma
      category $R/k$ is nonempty and connected. \qedhere
  \end{enumerate}
\end{definition}
\begin{proof}[Proof of equivalence]
  The equivalence characterizing an initial functor (or rather its dual) is
  proved in \cite[\mbox{Lemma 8.3.4}]{riehl2014}.
\end{proof}

The limit characterization of initiality is useful because one can often see at
a glance whether a given functor is initial just by imagining what will happen
to a limit of sets. The characterization based on connectedness is useful in
being concrete and combinatorial. To spell this out, given a functor
$R: \cat{J} \to \cat{K}$ and an object $k \in \cat{K}$, the \define{comma category}
$R/k$ has
\begin{itemize}
  \item as objects, an object $j \in \cat{J}$ together with a morphism
    $f: R(j) \to k$ in $\cat{K}$;
  \item as morphisms $(j, f) \to (j', f')$, morphisms $h: j \to j'$ in $\cat{J}$
    making the triangle commute:
    \begin{equation*}
      \begin{tikzcd}[column sep=small]
        {R(j)} && {R(j')} \\
        & k
        \arrow["{R(h)}", from=1-1, to=1-3]
        \arrow["f"', from=1-1, to=2-2]
        \arrow["{f'}", from=1-3, to=2-2]
      \end{tikzcd}.
    \end{equation*}
\end{itemize}
A category is \define{connected} if any two of its objects can be connected by a
zig-zag of morphisms, i.e., there is an \emph{undirected} path between them. In
the common case that the categories $\cat{J}$ and $\cat{K}$ are freely and
finitely generated, proving initiality---checking that each comma category $R/k$
is nonempty and connected---reduces to a combinatorial problem.

\begin{example}[Initial object as initial functor]
  Let $i$ be an object in a category $\cat{J}$. The corresponding functor
  $i: \cat{1} \to \cat{J}$ is initial if and only if for every $j \in \cat{J}$,
  the discrete category $i/j \cong \Hom_{\cat{J}}(i,j)$ is nonempty and
  connected, i.e., $i$ is an initial object in $\cat{J}$. This fact explains the
  term ``initial,'' which might otherwise appear backwards since initial
  functors are about limits, not colimits.
\end{example}

Initial functors are a well-behaved class of functors, closed under many
operations.
\begin{proposition}[Closure properties of initial functors]
  \label{prop:initial-functor-closure-properties}
  Let $F$ and $G$ be composable functors.
  \begin{enumerate}[(i),noitemsep]
    \item \emph{Composition}: If $F$ and $G$ are initial, then so is $G \circ F$.
    \item \emph{Cancellation}: If $F$ and $G \circ F$ are initial, then so is $G$.
    \item \emph{Equivalences}: Every full and essentially surjective
      functor---in particular, every equivalence of categories---is initial.
    \item \emph{Pushouts}: Initial functors are stable under pushout.
    \item \emph{Coproducts}: Initial functors are closed under coproducts (in
      $\cat{Cat}^{\to}$).
  \end{enumerate}
\end{proposition}

As the next subsection shows, nearly all of these properties follow from the
foundational theorem about initial functors: that initial functors form the left
class of morphisms in a factorization system on $\cat{Cat}$.

\subsection{Comprehensive factorization of functors}

The \emph{comprehensive factorization system}, discovered by Street and Walters
\cite{street1973}, is an orthogonal factorization system on the category of
small categories. Its left class of morphisms are the initial functors. Its
right class are the discrete opfibrations, whose definition we now recall.

\begin{definition}[Discrete opfibration]
  \label{def:discrete-opfibration}
  A functor $P: \cat{E} \to \cat{C}$ is a \define{discrete opfibration}, here
  denoted $P: \cat{E} \dopfto \cat{C}$, if for every morphism $f: X \to Y$ in
  $\cat{C}$ and every object $x \in \cat{E}$ with $P(x) = X$, there exists a
  unique morphism $\bar f: x \to y$ in $\cat{E}$ such that $P(\bar f) = f$:
  \begin{equation*}
    \begin{tikzcd}
      x & {} & {\operatorname{cod}(\bar f)} \\
      X & {} & Y
      \arrow["{\bar f}", dashed, from=1-1, to=1-3]
      \arrow["P"', between={0.2}{0.8}, maps to, from=1-2, to=2-2]
      \arrow["f"', from=2-1, to=2-3]
    \end{tikzcd}
  \end{equation*}
  Note that the codomain $y = \cod(\bar f)$ of $\bar f$ is part of the data
  required to exist uniquely.
\end{definition}

\begin{theorem}[Comprehensive factorization \cite[\mbox{Theorems 3 \& 4}]{street1973}]
  \label{thm:comprehensive-factorization}
  Initial functors and discrete opfibrations form the left and right classes,
  respectively, of an orthogonal factorization system on $\cat{Cat}$, the
  category of small categories and functors.
\end{theorem}

Using the description of an orthogonal factorization system in terms of unique
factorization up to unique isomorphism (see \cite[\mbox{Theorem 3.7}]{joyal2020}
or \cite[\mbox{Lemma 1.8}]{riehl2008}), the theorem says that any functor
$F: \cat{C} \to \cat{D}$ admits a factorization
\begin{equation*}
  \begin{tikzcd}
    {\cat{C}} && {\cat{D}} \\
    & {\cat{E}}
    \arrow["F", from=1-1, to=1-3]
    \arrow["E"', two heads, from=1-1, to=2-2]
    \arrow["P"', tail, from=2-2, to=1-3]
  \end{tikzcd}
\end{equation*}
as an initial functor followed by a discrete opfibration, and such a
factorization is unique in the sense that any other factorization
$(\cat{E}',E',P')$ is related by a unique isomorphism making both triangles
commute:
\begin{equation*}
  \begin{tikzcd}[row sep=small]
    & {\cat{E}} & \\
    {\cat{C}} && {\cat{D}} \\
    & {\cat{E}'}
    \arrow["P", tail, from=1-2, to=2-3]
    \arrow["{\exists!}"', from=1-2, to=3-2]
    \arrow["\cong"{marking, allow upside down}, shift left=2, draw=none, from=1-2, to=3-2]
    \arrow["E", two heads, from=2-1, to=1-2]
    \arrow["{E'}"', two heads, from=2-1, to=3-2]
    \arrow["{P'}"', tail, from=3-2, to=2-3]
  \end{tikzcd}.
\end{equation*}

Later, we will explain how to actually construct the comprehensive
factorization. Its existence already has many consequences, starting with:

\begin{proof}[Proof of \cref{prop:initial-functor-closure-properties}]
  Nearly all of the stated closure properties of initial functors follow from
  closure properties of the left class in an orthogonal factorization system,
  applied to the comprehensive factorization system
  (\cref{thm:comprehensive-factorization}). Specifically, the cancellation
  property is \cite[\mbox{Proposition 2.2}]{joyal2020} or \cite[\mbox{Lemma
    1.13}]{riehl2008}, and closure under colimits is \cite[\mbox{Theorem
    1.17}]{riehl2008}. The only property not covered by generalities about
  factorization systems is that full and essentially surjective functors are
  initial, proved in \cite[\mbox{Lemma 9.8}]{patterson2023}.
\end{proof}

The comprehensive factorization extends from categories to diagrams in a fixed
category.

\begin{corollary}[Comprehensive factorization for diagrams]
  \label{cor:comprehensive-factorization-diagrams}
  Let $\cat{C}$ be any category. The slice category $\cat{Cat}/\cat{C}$, or
  \emph{category of diagrams} in $\cat{C}$, has an orthogonal factorization
  system whose left and right classes are the morphisms underlain by initial
  functors and discrete opfibrations, respectively.
\end{corollary}
\begin{proof}
  This follows from \cref{thm:comprehensive-factorization} and the fact that
  orthogonal factorization systems pass to slices (and coslices)
  \cite[\S{1}]{joyal2020}.
\end{proof}

Initial functors can also serve as weak equivalences in the sense of homotopy
theory. Though not used directly in the main text, the result is worth recording
since it too follows from the comprehensive factorization.

\begin{corollary}[Calculi of fractions]
  \label{cor:calculus-of-fractions}
  As a class of weak equivalences in $\cat{Cat}$, the initial functors admit a
  \emph{calculus of left fractions} in the sense of Gabriel--Zisman
  \cite[\S{I.2}]{gabriel1967}.

  For any category $\cat{C}$, the same is true of the morphisms in
  $\cat{Cat}/\cat{C}$ underlain by initial functors.
\end{corollary}
\begin{proof}
  The two statements follow from \cref{thm:comprehensive-factorization} and
  \cref{cor:comprehensive-factorization-diagrams}, combined with the fact that
  the left class in an orthogonal factorization system admits a calculus of
  left fractions \cite[\mbox{Example 2.4}]{kurz2005}.
\end{proof}

That initial morphisms between diagrams admit a calculus of left fractions was
first proved in \cite[\S{3}, p.\ 32]{pare1973}. The main idea behind that paper
is the subject of the next subsection.

\subsection{The image copresheaf of a diagram}

The two seminal early works on initial functors, both published in 1973, are by
Street and Walters \cite{street1973}, on the comprehensive factorization of a
functor, and by Par{\'e} \cite{pare1973}, on the image copresheaf generated by a
diagram.\footnote{The name ``image copresheaf'' is borrowed from the more recent
  exposition \cite[\S{3}]{perrone2022}.} Their connection, reviewed in this
subsection, lies in how the image copresheaf is used to construct the
comprehensive factorization, analogously to how the image of a function is used
to construct the epi-mono factorization.

Many different diagrams in a category can have the same limit. The image
copresheaf captures all of the functorial limit properties of a diagram while
forgetting the accidental properties of its shape \cite[cf.][\S{0}]{pare1973}.

\begin{definition}[Image copresheaf]
  \label{def:image-copresheaf}
  The \define{image copresheaf} of a diagram $D: \cat{J} \to \cat{C}$ is the
  copresheaf on $\cat{C}$ defined by
  \begin{equation*}
    \Im D \coloneqq \Big(\cat{C} \xto{D/(-)} \cat{Cat} \xto{\pi_0} \Set\Big),
  \end{equation*}
  where $D/(-)$ is the comma category construction and $\pi_0$ is the connected
  components functor.
\end{definition}

So, on objects, $\Im(D)(c) \coloneqq \pi_0(D/c)$ for each $c \in \cat{C}$. The
similarity with the comma characterization of an initial functor
(\cref{def:initial-functor}) is striking. To be more precise, a functor
$R: \cat{J} \to \cat{K}$ is initial if and only if $\Im R$ is naturally isomorphic
to $1$, the terminal copresheaf on $\cat{K}$.

Like the notion of initial functor, the image copresheaf admits several
characterizations:

\begin{theorem}[Guises of image copresheaf]
  \label{thm:image-copresheaf}
  Let $D: \cat{J} \to \cat{C}$ be a diagram. Then the following copresheaves on
  $\cat{C}$ are naturally isomorphic:
  \begin{enumerate}[(i)]
    \item $\Im D \coloneqq \pi_0(D/-)$;
    \item the colimit of representables $\colim_{\cat{J}^\op} (y \circ D^\op)$, where
      $y: \cat{C}^\op \to \Set^{\cat{C}}$ is the (contravariant) Yoneda embedding
      of $\cat{C}$;
    \item the left Kan extension $\Lan_D(1)$, where $1$ is the terminal copresheaf on $\cat{J}$;
    \item a representing object for the representable functor
      \begin{equation*}
        \textstyle \lim_{\cat{J}} (- \circ D):
          \Set^{\cat{C}} \to \Set,\ F \mapsto \lim_{\cat{J}} (F \circ D).
      \end{equation*}
    \end{enumerate}
\end{theorem}
\begin{proof}
  The isomorphisms among (i)--(iii) are (the dual of) \cite[Theorem 2.1 and
  Corollary 2.8]{pare1973}.

  We prove that (iii) and (iv) are isomorphic. Given a functor
  $F: \cat{C} \to \Set$, we have
  \begin{equation*} \textstyle
    \lim_{\cat{J}} (F \circ D) \cong \Set^{\cat{J}}(1, F \circ D)
      \cong \Set^{\cat{C}}(\Lan_D(1), F).
  \end{equation*}
  Hence, the functor $\lim_{\cat{J}} (- \circ D): \Set^{\cat{C}} \to \Set$ is
  represented by $\Lan_D(1) \cong \Im D$ (which, by the Yoneda lemma, is unique
  up to isomorphism).
\end{proof}

The fundamental theorem about the image copresheaf is that it classifies
diagrams with the same functorial limit properties, which is reinterpreted in
\cref{sec:soundness-completeness} as a soundness and completeness result for
equational reasoning with diagrams.

\begin{theorem}[Isomorphism of image copresheaves]
  \label{thm:image-copresheaf-isomorphism}
  Let $D: \cat{J} \to \cat{C}$ and $D': \cat{J}' \to \cat{C}$ be diagrams in a
  category $\cat{C}$. The following are equivalent:
  \begin{enumerate}[(i)]
    \item the copresheaves $\Im D$ and $\Im D'$ are naturally isomorphic;
    \item the diagrams $D$ and $D'$ are weakly equivalent
      (\cref{sec:weak-equivalence}), i.e., they can be connected by a zig-zag in
      $\cat{Cat}/\cat{C}$ with each morphism underlain by an initial functor;
    \item for every category $\cat{D}$ and functor $F: \cat{C} \to \cat{D}$, the
      limits $\lim_{\cat{J}} (F \circ D)$ and $\lim_{\cat{J}'} (F \circ D')$ are
      isomorphic, each existing if the other does;
    \item for every set-valued functor $F: \cat{C} \to \Set$, the limits
      $\lim_{\cat{J}} (F \circ D)$ and $\lim_{\cat{J}'} (F \circ D')$ are in
      bijection, naturally in $F$.
  \end{enumerate}
\end{theorem}
\begin{proof}
  We have been careless about size issues, but in fact the theorem is true
  regardless of whether $\cat{C}$ and $\cat{D}$ are required to be \emph{small}
  categories or range over all \emph{locally small} categories. In the case of
  small categories, the equivalence of (i) and (iii) is \cite[\mbox{Theorem
    3.2}]{pare1973}, while the equivalence with (ii) is implicit in \cite[end of
  \S{3}]{pare1973}. In the case of locally small categories, the equivalence of
  (i)--(iii) appears verbatim as \cite[\mbox{Proposition 3.9}]{perrone2022}.

  As for the equivalence of (i) and (iv), we have $\Im D \cong \Im D'$ if and
  only if $\lim_{\cat{J}} (- \circ D) \cong \lim_{\cat{J}'} (- \circ D')$, by
  \cref{thm:image-copresheaf} and the Yoneda lemma.
\end{proof}

Finally, we recall how the image copresheaf appears in the construction of the
comprehensive factorization of a functor.

\begin{construction}[Comprehensive factorization \cite{street1973}]
  \label{con:comprehensive-factorization}
  Let $F: \cat{C} \to \cat{D}$ be a functor out of a small category $\cat{C}$. The
  image copresheaf of $F$ is the copresheaf $\Im F \cong \Lan_F 1$ on $\cat{D}$, the
  left Kan extension of the terminal copresheaf on $\cat{C}$ along $F$. By the
  universal property of its category of elements as a comma category,
  ${\El(\Im F) \cong 1/\Im F}$, there exists a unique functor
  ${E: \cat{C} \to \El(\Im F)}$ factorizing $F$ as
  \begin{equation*}
    \begin{tikzcd}[column sep=small]
      {\cat{C}} && {\cat{D}} \\
      & {\El(\Im F)}
      \arrow["F", from=1-1, to=1-3]
      \arrow["E"', from=1-1, to=2-2]
      \arrow["{\pi_{\Im F}}"', from=2-2, to=1-3]
    \end{tikzcd}
  \end{equation*}
  and equating the 2-cells:
  \begin{equation*}
    \begin{tikzcd}
      {\cat{C}} && \\
      & {\El(\Im F)} & {\cat{D}} \\
      & {\cat{1}} & \Set
      \arrow["E", from=1-1, to=2-2]
      \arrow["{\pi_{\Im F}}", from=2-2, to=2-3]
      \arrow["{!}"', from=2-2, to=3-2]
      \arrow["{\Im F}", from=2-3, to=3-3]
      \arrow[between={0.2}{0.8}, Rightarrow, nfold, from=3-2, to=2-3]
      \arrow["1"', from=3-2, to=3-3]
    \end{tikzcd}
    \quad=\quad
    \begin{tikzcd}
      {\cat{C}} & {\cat{D}} \\
      {\cat{1}} & \Set
      \arrow["F", from=1-1, to=1-2]
      \arrow["{!}"', from=1-1, to=2-1]
      \arrow["{\Im F\,\cong\,\Lan_F 1}", from=1-2, to=2-2]
      \arrow["\eta", between={0.2}{0.8}, Rightarrow, nfold, from=2-1, to=1-2]
      \arrow["1"', from=2-1, to=2-2]
    \end{tikzcd}.
  \end{equation*}
  In this factorization, $\pi_{\Im F}$ is a discrete opfibration, and Street and Walters
  show that $E$ is initial. So ${F = \pi_{\Im F} \circ E}$ is the \define{comprehensive
    factorization} of $F$.
\end{construction}

\section{Cartesian initial functors}
\label{sec:cartesian-initial-functors}

In \cref{sec:cartesian-setting}, we consider diagrams, limits, and initial
functors in the setting of categories with finite products. This appendix shows
that the well-established theory of these concepts for ordinary categories,
reviewed in \cref{sec:initial-functors}, extends smoothly to cartesian
categories. Although we cannot locate a few of the statements directly in the
literature, this material should be considered essentially known.

\subsection{Initial functors between cartesian categories}

For initial functors between \emph{cartesian} categories, the defining
characterizations involving limits (\cref{def:initial-functor}) still hold when
confined to cartesian diagrams:

\begin{proposition}
  \label{prop:cartesian-initial-functor}
  Let $R: \cat{J} \to \cat{K}$ be a functor between cartesian categories.
  The following are equivalent:
  \begin{enumerate}[(i)]
    \item $R$ is initial;
    \item for every cartesian diagram $D: \cat{K} \to \cat{C}$ in a cartesian category
      $\cat{C}$, the canonical comparison
      $\lim_{\cat{K}} D \to \lim_{\cat{J}} (D \circ R)$ is an isomorphism, each limit
      existing if the other does;
    \item for every cartesian diagram $D: \cat{K} \to \Set$, the canonical comparison
      $\lim_{\cat{K}} D \to \lim_{\cat{J}} (D \circ R)$ is a bijection.
  \end{enumerate}
\end{proposition}
\begin{proof}
  That (i) implies (ii) and (ii) implies (iii) is immediate. We prove that (iii)
  implies (i). The equivalence of (i) and (ii) can also be deduced, in dual
  form, from the ``relative theory'' in \cite[\S\S{5--6}]{pare1973}.

  Let $[\cat{K}, \Set]_{\mathrm{cart}}$ denote the full subcategory of
  $\Set^{\cat{K}}$ spanned by the cartesian functors. Then condition (iii)
  furnishes a natural isomorphism
  $\lim_{\cat{K}} (-) \cong \lim_{\cat{J}}(- \circ R)$ between functors
  $[\cat{K}, \Set]_{\mathrm{cart}} \to \Set$. Now, the terminal copresheaf on
  $\cat{K}$ is trivially cartesian, and the image copresheaf $\Im R$ is
  cartesian by \cref{cor:cartesian-image-copresheaf}, proved below. Thus,
  $\Im R$ and $1$ are representing objects for $\lim_{\cat{J}}(- \circ R)$ and
  $\lim_{\cat{K}}(-)$, respectively, using the representation in
  \cref{thm:image-copresheaf}. Apply the Yoneda lemma to conclude that
  $\Im R \cong 1$ and hence that $R$ is initial.
\end{proof}

Thus, in defining a \define{cartesian initial functor} to be simply a functor
that is both cartesian and initial, we need not worry about which adjective
comes first.

\subsection{Comprehensive factorization of cartesian functors}

The foundational theorem about initial functors is the existence of the
comprehensive factorization system \cite{street1973}, an orthogonal
factorization system on $\cat{Cat}$ having initial functors and discrete
opfibrations as its left and right classes
(\cref{thm:comprehensive-factorization}). Perhaps surprisingly, the
comprehensive factorization works just as well for cartesian functors as it does
for ordinary functors. The precise statement is:

\begin{theorem}[Cartesian comprehensive factorization]
  \label{thm:cartesian-comprehensive-factorization}
  The comprehensive factorization system on $\cat{Cat}$ restricts to an
  orthogonal factorization system on $\cat{Cart}$, the category of small
  categories with finite products and finite-product-preserving functors.

  Moreover, this factorization system on $\cat{Cart}$ lifts to an orthogonal
  factorization system on $\cat{Cart}_{\mathrm{str}}$, the category of small
  categories with chosen finite products and functors strictly preserving the
  chosen products.
\end{theorem}

We do not know of a source for this theorem in the literature, though we suspect
that it is folklore since it follows straightforwardly from other,
well-documented facts. In any case, the purpose of this subsection is to present
a self-contained proof. To do so, we show that the explicit construction of the
comprehensive factorization (\cref{con:comprehensive-factorization}), when
applied to a cartesian functor, is a factorization in which the intermediate
category and both functors are automatically cartesian, without further
assumptions.

The analysis of the right class---discrete opfibrations---is simpler, so we
start there. Recall that discrete opfibrations correspond to copresheaves via
the elements construction \cite[\S{2.4}]{riehl2016}. A simple result says that
the elements construction preserves cartesianness.

\begin{lemma}[Products in elements construction]
  \label{lem:cartesian-elements-construction}
  Let $X: \cat{C} \to \Set$ be a cartesian set-valued functor on a cartesian
  category $\cat{C}$. Then the category of elements, $\El(X)$, is cartesian, and
  its canonical projection, the discrete opfibration ${\pi_X: \El(X) \to \cat{C}}$,
  both preserves and reflects finite products.
\end{lemma}
\begin{proof}
  The product of objects $(c_i, x_i)_{i \in I}$ in $\El(X)$ is given, in its
  first component, by the product in $\cat{C}$, and, in its second component, by
  pairing:
  \begin{equation*}
    \left(
      \prod_{i \in I} c_i,\
      \langle x_i \rangle_{i \in I} \in \prod_{i \in I} X(c_i)
        \cong X\big({\textstyle\prod_{i \in I} c_i}\big)
    \right).
  \end{equation*}
  Product projections in $\El(X)$ are simply the product projections in
  $\cat{C}$. The canonical projection $\pi_X: {(c, x) \mapsto c},\, {f \mapsto f}$ clearly
  preserves finite products, and it also reflects them since a tuple is uniquely
  determined by its components (or, apply \cite[\mbox{Lemma
    6.8}]{patterson2023}).
\end{proof}

Since the comprehensive factorization of a functor involves a left Kan
extension, we ask how left Kan extensions interact with finite products. The
crucial concept in answering this question is \emph{siftedness}
\cite{gabriel1971}---or rather its dual, \emph{cosiftedness}, characterized by
the following definition/theorem.

\begin{definition}[Cosiftedness]
  \label{def:cosifted}
  A small category $\cat{J}$ is \define{cosifted} if any of the following
  equivalent conditions holds:
  \begin{enumerate}[(i),noitemsep]
    \item in $\Set$, colimits of shape $\cat{J}^\op$ commute with finite
      products;
    \item $\cat{J}$ is nonempty and the diagonal functor
      $\Delta_{\cat{J}}: \cat{J} \to \cat{J} \times \cat{J}$ is initial;
    \item $\cat{J}$ is nonempty and for any objects $j_1, j_2 \in \cat{J}$, the
      category $\cat{Span}_{\cat{J}}(j_1, j_2)$ of spans in $\cat{J}$ with feet
      $j_1$ and $j_2$ is nonempty and connected.
  \end{enumerate}
  A \define{cosifted colimit} is a colimit whose shape is a cosifted category.
\end{definition}
\begin{proof}[Proof of equivalence]
  The (dual of the) equivalence between (i) and (ii) originates in
  \cite{gabriel1971}; see also \cite[\mbox{Theorem 2.15}]{adamek2011}. The
  equivalence between (ii) and (iii) follows from the defining characterization
  of an initial functor (\cref{def:initial-functor}).
\end{proof}

\begin{lemma}
  \label{lem:cartesian-category-cosifted}
  A small category with finite products is cosifted.
\end{lemma}
\begin{proof}
  Suppose $\cat{J}$ has finite products. Then $\cat{J}$ is nonempty, having a
  terminal object, and for any objects $j_1, j_2 \in \cat{J}$, the category
  $\cat{Span}_{\cat{J}}(j_1, j_2)$ has a terminal object---the product cone
  ${j_1 \from j_1 \times j_2 \to j_2}$, by the universal property of the
  product---so is in particular nonempty and connected.
\end{proof}

The surprising fact is that the left Kan extension of a cartesian set-valued
functor along \emph{any} functor (not necessarily cartesian) is again cartesian.
This result seems to first appear in the literature as \cite[\mbox{Theorem 1.5
  and Example 3.1}]{borceux1977}, in the more general setting of enriched
categories. We give a direct and elementary proof.

\begin{theorem}[Product-preserving Kan extensions]
  \label{thm:product-preserving-kan-extensions}
  Let $\cat{C}$ be a small cartesian category, and let $X: \cat{C} \to \Set$ be a
  cartesian set-valued functor. For any functor $F: \cat{C} \to \cat{D}$ into a
  cartesian category $\cat{D}$, the left Kan extension
  $\Lan_F X: \cat{D} \to \Set$ is cartesian.
\end{theorem}
\begin{proof}
  Decomposing $X$ as a canonical colimit of representables, the left Kan
  extension of $X$ has the formula:
  \begin{equation*}
    \Lan_F X \cong \colim_{(c,x) \in {\El(X)}^\op} \cat{D}(F(c), -).
  \end{equation*}
  Because the representables $\cat{D}(F(c), -): \cat{D} \to \Set$ preserve
  finite products, $\Lan_F X$ will preserve finite products when $\El(X)$ is
  cosifted. But that is the case by
  \cref{lem:cartesian-elements-construction,lem:cartesian-category-cosifted}.
\end{proof}

The following corollary also appears, in dual and slightly weaker form, as
\cite[\mbox{Proposition 6.5}]{pare1973}.

\begin{corollary}
  \label{cor:cartesian-image-copresheaf}
  Let $D: \cat{J} \to \cat{C}$ be a diagram whose shape and target are both
  cartesian categories. Then the image copresheaf of $D$ is cartesian.
\end{corollary}
\begin{proof}
  Since $1 \times 1 \cong 1$, the terminal copresheaf on $\cat{J}$ is cartesian.
  Thus, by \cref{thm:product-preserving-kan-extensions}, so is
  $\Im D \cong \Lan_D 1$.
\end{proof}

We now have all the ingredients needed to prove that the comprehensive
factorization on $\cat{Cat}$ restricts to $\cat{Cart}$, the category of small
cartesian categories and cartesian functors.

\begin{proof}[Proof of \cref{thm:cartesian-comprehensive-factorization}]
  Let $F: \cat{C} \to \cat{D}$ be a functor out of a small category $\cat{C}$,
  and adopt the notation of \cref{con:comprehensive-factorization}. If $\cat{C}$
  and $\cat{D}$ are cartesian, then so is $\Im F$ by
  \cref{cor:cartesian-image-copresheaf}; hence, so are $\El(\Im F)$ and $\pi_{\Im F}$ by
  \cref{lem:cartesian-elements-construction}. Moreover, if $F$ is itself
  cartesian, then so is $E$, since the composite ${\pi_{\Im F} \circ E = F}$ preserves
  finite products and, by \cref{lem:cartesian-elements-construction}, $\pi_{\Im F}$
  reflects them.

  Thus, when a functor $F: \cat{C} \to \cat{D}$ belongs to $\cat{Cart}$, the
  factorization ${F = \pi_{\Im F} \circ E}$ in $\cat{Cat}$ restricts to $\cat{Cart}$.
  Any morphism between factorizations also preserves finite products, since
  cartesian discrete opfibrations reflect them \cite[\mbox{Lemma
    6.8}]{patterson2023}. This proves the first part of the theorem, using the
  characterization of a factorization system in terms of unique factorization up
  to unique isomorphism \cite[\mbox{Theorem 3.7}]{joyal2020}.

  Now suppose that $F: \cat{C} \to \cat{D}$ belongs to
  $\cat{Cart}_{\mathrm{str}}$, the category of small categories with chosen
  finite products. There is a \emph{unique} choice of finite products in the
  elements construction $\El(\Im F)$ making the projection
  $\pi_{\Im F}: \El(\Im F) \to \cat{D}$ strictly preserve products. Since
  $F = \pi_{\Im F} \circ E$ also strictly preserves products, $E$ must do so
  too, as its product comparisons lie over identities in $\cat{D}$, hence must
  be identities by the discrete opfibration property. Similarly, any morphism
  between such factorizations must strictly preserve products. Thus the
  orthogonal factorization system lifts through the forgetful functor
  $\cat{Cart}_{\mathrm{str}} \to \cat{Cart}$.
\end{proof}

\subsection{Limits of cartesian diagrams}

In \cref{sec:cartesian-setting}, we take limits of diagrams whose shapes
themselves have finite products, a not-too-common situation in category theory.
This subsection shows that such limits are innocuous provided the diagram itself
preserves finite products.

We say that a diagram is \define{cartesian} if its shape and target have finite
products and if it, as a functor, preserves finite products. A first observation
about cartesian diagrams is:

\begin{lemma}[Cones over a cartesian diagram]
  \label{lem:cones-cartesian-diagram}
  A cone over a cartesian diagram is uniquely determined at a product of objects
  by its legs at each factor.
\end{lemma}
\begin{proof}
  Let $\lambda: c \To D$ be a cone over a cartesian diagram $D: \cat{J} \to \cat{C}$.
  For any finite list of objects $x_1, \dots, x_k$ in $\cat{J}$, pairing the
  naturality triangles at each projection $\pi_i: x_1 \times \cdots \times x_k \to x_i$ gives a
  commuting triangle
  \begin{equation*}
    \begin{tikzcd}[column sep=large]
      c & {D(x_1 \times \cdots \times x_k)} \\
      & {D(x_1) \times \cdots \times D(x_k)}
      \arrow["{\lambda_{x_1 \times \cdots \times x_k}}", from=1-1, to=1-2]
      \arrow["{\langle \lambda_{x_i} \rangle_{i=1}^k}"', from=1-1, to=2-2]
      \arrow["{\langle D(\pi_i) \rangle_{i=1}^k}", from=1-2, to=2-2]
      \arrow["\cong"{marking, allow upside down}, shift left=2, draw=none, from=2-2, to=1-2]
    \end{tikzcd},
  \end{equation*}
  in which the comparison map is an isomorphism since $D$ preserves finite
  products.
\end{proof}

The next two propositions show that cartesian diagrams freely generated by two
different means have the expected cones, hence the expected limits.

\begin{proposition}[Cones over a free cartesian diagram]
  \label{prop:cones-free-cartesian-diagram}
  Let $D: \cat{J} \to \cat{C}$ be a free cartesian diagram, with shape
  $\cat{J} = F(M)$ freely generated by a multigraph $M$. Then any family of
  morphisms ${\lambda_x: c \to D(x)}$ in $\cat{C}$ indexed by vertices $x$ of
  $M$, making the square commute
  \begin{equation*}
    \begin{tikzcd}[column sep=large]
      c & {D(x_1) \times \cdots \times D(x_k)} \\
      {D(y)} & {D(x_1 \times \cdots \times x_k)}
      \arrow["{\langle \lambda_{x_1}, \dots, \lambda_{x_k} \rangle}", from=1-1, to=1-2]
      \arrow["{\lambda_y}"', from=1-1, to=2-1]
      \arrow["\cong", from=1-2, to=2-2]
      \arrow["{D(f)}", from=2-2, to=2-1]
    \end{tikzcd}
  \end{equation*}
  for every multiedge $f: (x_1, \dots, x_k) \to y$ in $M$, extends uniquely to a
  cone $\lambda: c \To D$ over $D$.

  Moreover, this bijective correspondence with cones over $D$ is natural in the
  apex $c$.
\end{proposition}

For the next proposition, let $\ProdComp(\cat{J})$ denote the free completion of a
category $\cat{J}$ under finite products, given explicitly by the contravariant
family construction \cite[e.g.,][\mbox{Remark 1.15}]{adamek2011}.

\begin{proposition}[Limit indexed by a free product completion]
  \label{prop:limit-free-product-completion}
  Let $D: \cat{J} \to \cat{C}$ be a diagram in a cartesian category $\cat{C}$, and
  let $\hat D$ be its essentially unique extension to a cartesian diagram:
  \begin{equation*}
    \begin{tikzcd}
      & {\ProdComp(\cat{J})} \\
      {\cat{J}} & {\cat{C}}
      \arrow["{\hat D}", dashed, from=1-2, to=2-2]
      \arrow["{\eta_{\cat{J}}}", from=2-1, to=1-2]
      \arrow["D"', from=2-1, to=2-2]
    \end{tikzcd}.
  \end{equation*}
  Then cones over $D$ and $\hat D$ are in bijection, naturally in the apex; in
  particular, $\lim D \cong \lim \hat D$, each limit existing if the other does.
\end{proposition}

Both propositions can be proved by inspecting explicit constructions of the free
completions. We take a more abstract approach that uses their universal
properties instead.

To do so, we must change our perspective on cones, viewing them not as natural
transformations but as functors solving a lifting problem. Fix a diagram
$D: \cat{J} \to \cat{C}$. For any object $c$ in $\cat{C}$, let
$\cod: c/\cat{C} \to \cat{C}$ be the codomain functor out of the coslice category
under $c$. Then cones with apex $c$ over $D$ correspond exactly with lifts of
$D$ through the codomain functor:
\begin{equation} \label{eq:cones-and-lifts}
  \lambda: c \To D
  \qquad\leftrightsquigarrow\qquad
  \begin{tikzcd}
    & {c/\cat{C}} \\
    {\cat{J}} & {\cat{C}}
    \arrow["\cod", from=1-2, to=2-2]
    \arrow["{\bar D}", dashed, from=2-1, to=1-2]
    \arrow["D"', from=2-1, to=2-2]
  \end{tikzcd}.
\end{equation}
This perspective is taken consistently in \cite[esp.\ \S{2}]{patterson2023}.

As preparation, we restate and reprove the first lemma \cite[cf.][\mbox{Lemma
  6.7}]{patterson2023}.

\begin{lemma}[{\cref{lem:cones-cartesian-diagram}, for lifts}]
  \label{lem:lifts-cartesian-diagram}
  When $D: \cat{J} \to \cat{C}$ is a cartesian diagram, the lifting problem
  \eqref{eq:cones-and-lifts} is the same whether it is taken in the category
  $\cat{Cat}$ or $\cat{Cart}$.
\end{lemma}
\begin{proof}
  Let $c$ be an object of $\cat{C}$. Since $\cat{C}$ is cartesian, so is the
  representable $\cat{C}(c,-): \cat{C} \to \Set$; hence, by
  \cref{lem:cartesian-elements-construction}, the canonical projection
  $\pi: \El(\cat{C}(c,-)) \to \cat{C}$---which is simply
  $\cod: c/\cat{C} \to \cat{C}$---reflects finite products. Since $D$ preserves
  finite products, it follows that any lift $\bar D$ of $D$ through $\cod$ (in
  the category $\cat{Cat}$) also preserves finite products.
\end{proof}

To see the idea of the proof of \cref{prop:cones-free-cartesian-diagram}, imagine
there were a free-forgetful adjunction
$F: \MGraph \rightleftarrows \cat{Cart}: U$ between multigraphs and cartesian
categories. Let the cartesian diagram $D: F(M) \to \cat{C}$ and the multigraph
morphism $d: M \to U(\cat{C})$ correspond under the adjunction. Then naturality
of the hom-set bijections \cite[\mbox{Lemma 4.1.3}]{riehl2016} implies that
lifts of $D$ and lifts of $d$ correspond bijectively:
\begin{equation*}
  \begin{tikzcd}
    & {c/\cat{C}} \\
    {F(M)} & {\cat{C}}
    \arrow["\cod", from=1-2, to=2-2]
    \arrow["{\bar D}", dashed, from=2-1, to=1-2]
    \arrow["D"', from=2-1, to=2-2]
  \end{tikzcd}
  \qquad\leftrightsquigarrow\qquad
  \begin{tikzcd}
    & {U(c/\cat{C})} \\
    M & {U(\cat{C})}
    \arrow["{U(\cod)}", from=1-2, to=2-2]
    \arrow["{\bar d}", dashed, from=2-1, to=1-2]
    \arrow["d"', from=2-1, to=2-2]
  \end{tikzcd}.
\end{equation*}
But lifts of $D$ correspond to cones over $D$ (noting that lifts in $\cat{Cat}$
and $\cat{Cart}$ are the same, by \cref{lem:lifts-cartesian-diagram}), and lifts
of $d$ correspond to the $M$-indexed families in the proposition statement. So the
first part of \cref{prop:cones-free-cartesian-diagram} would be proved.

The problem with this argument is that there is not, in fact, a (1-categorical)
adjunction between multigraphs and cartesian categories. The actual proof uses
the 2-categorical universal property discovered by Walters \cite{walters1989},
interesting in its own right.

\begin{proof}[Proof of \cref{prop:cones-free-cartesian-diagram}]
  To state a 2-categorical universal property for the free cartesian category,
  Walters \cite{walters1989} constructs the 2-category $\MGraph(\Cat)$ of
  multigraphs internal to $\Cat$ and a forgetful 2-functor
  $U: \Cart \to \MGraph(\Cat)$. Viewing an ordinary multigraph $M$ as a discrete
  object in $\MGraph(\Cat)$, there is a multigraph morphism $\eta_M: M \to U(F(M))$,
  pre-composition with which induces an equivalence of hom-categories:
  \begin{equation*}
    U(-) \circ \eta_M: \Cart(F(M), \cat{C})
      \xto{\sim} \MGraph(\Cat)(M, U(\cat{C})).
  \end{equation*}
  This is the better part of a biadjunction. In any case, for each object $c$ of
  $\cat{C}$, we obtain a commuting square
  \begin{equation*}
    \begin{tikzcd}[column sep=large]
      {\Cart(F(M),c/\cat{C})} & {\MGraph(\Cat)(M, U(c/\cat{C}))} \\
      {\Cart(F(M),\cat{C})} & {\MGraph(\Cat)(M, U(\cat{C}))}
      \arrow["{U(-) \circ \eta_M}", from=1-1, to=1-2]
      \arrow["{\cod \circ (-)}"', from=1-1, to=2-1]
      \arrow["{U(\cod) \circ (-)}", from=1-2, to=2-2]
      \arrow["{U(-) \circ \eta_M}"', from=2-1, to=2-2]
    \end{tikzcd}
  \end{equation*}
  in which the top and bottom functors are equivalences.

  Since $\cod$ is a discrete opfibration, post-composing with it gives a functor
  $\cod \circ (-): \Cat(\cat{J}, c/\cat{C}) \to \Cat(\cat{J}, \cat{C})$ that is itself
  a discrete opfibration, for any category $\cat{J}$. Moreover, since $\cod$ is
  cartesian, the post-composition functor
  $\cod \circ (-): \Cart(\cat{J}, c/\cat{C}) \to \Cart(\cat{J}, \cat{C})$ is also a
  discrete opfibration, for any cartesian category $\cat{J}$. These two facts
  are special cases of \cite[Theorems 5.2 and 6.6]{patterson2023}, respectively.
  Post-composition with $U(\cod)$ can also be shown to be a discrete
  opfibration.

  So, in the commuting square, the left and right functors are discrete
  opfibrations and the top and bottom functors are equivalences. In this
  situation, the top functor restricts to a bijection on fibers; here, that
  means that if a cartesian diagram $D: F(M) \to \cat{C}$ and a multigraph
  morphism $d: M \to U(\cat{C})$ correspond under the bottom equivalence, then the
  top equivalence puts their lifts, hence also their cones, in bijection.
  Moreover, the bijection is natural with respect to changes of total space
  $h^*: c/\cat{C} \to c'/\cat{C}$ given by pre-composing with a morphism
  $h: c' \to c$ in $\cat{C}$.
\end{proof}

\begin{proof}[Proof of \cref{prop:limit-free-product-completion}]
  The proof is similar to that of \cref{prop:cones-free-cartesian-diagram}.

  Let $\ProdComp: \Cat \rightleftarrows \Cart: U$ be the free-forgetful
  biadjunction between categories and categories with finite products. For any
  category $\cat{J}$ and cartesian category $\cat{C}$, pre-composition with the
  unit of the biadjunction, $\eta_{\cat{J}}: \cat{J} \to U(\ProdComp(\cat{J}))$,
  gives an equivalence of hom-categories
  \begin{equation*}
    U(-) \circ \eta_{\cat{J}}: \Cart(\ProdComp(\cat{J}), \cat{C})
      \xto{\sim} \Cat(\cat{J}, U(\cat{C})).
  \end{equation*}
  Then, for each object $c$ of $\cat{C}$, there is a commuting square of
  functors
  \begin{equation*}
    \begin{tikzcd}[column sep=large]
      {\Cart(\ProdComp(\cat{J}),c/\cat{C})} & {\Cat(\cat{J}, U(c/\cat{C}))} \\
      {\Cart(\ProdComp(\cat{J}),\cat{C})} & {\Cat(\cat{J}, U(\cat{C}))}
      \arrow["{U(-) \circ \eta_{\cat{J}}}", from=1-1, to=1-2]
      \arrow["{\cod \circ (-)}"', from=1-1, to=2-1]
      \arrow["{U(\cod) \circ (-)}", from=1-2, to=2-2]
      \arrow["{U(-) \circ \eta_{\cat{J}}}"', from=2-1, to=2-2]
    \end{tikzcd}
  \end{equation*}
  in which the top and bottom are equivalences and the left and right are
  discrete opfibrations \cite[Theorems 5.2 and 6.6]{patterson2023}. The rest of
  the proof follows the previous one.
\end{proof}

\end{document}